\documentclass[10pt]{amsart}
\usepackage[utf8]{inputenc}
\usepackage[T1]{fontenc}
\usepackage{indentfirst}
\usepackage{textcomp}
\usepackage[english]{babel}
\usepackage{amsmath,amsfonts,amssymb,amsopn,amscd,amsthm}
\usepackage{mathbbol}
\usepackage{mathrsfs}
\usepackage{stmaryrd}
\usepackage{graphicx}
\usepackage[dvipsnames]{xcolor}
\usepackage{pstricks,pst-plot}
\usepackage{xy}  \xyoption{all}
\usepackage{pb-diagram,pb-xy}
\usepackage{enumerate}
\usepackage[colorlinks=true,citecolor=DarkOrchid,linkcolor=NavyBlue]{hyperref}

\newtheorem{definitionproposition}{Definition - Proposition}
\newtheorem{theorem}{Theorem}
\newtheorem{proposition}[theorem]{Proposition}
\newtheorem{corollary}[theorem]{Corollary}
\newtheorem{lemma}[theorem]{Lemma}
\theoremstyle{remark}
\newtheorem{example}{Example}
\newtheorem{examples}[example]{Examples}
\newtheorem{remark}{Remark}

\newcommand{\E}{\mathrm{e}}
\newcommand{\N}{\mathbb{N}}

\newcommand{\R}{\mathbb{R}}
\newcommand{\C}{\mathbb{C}}
\newcommand{\Func}{\mathcal{F}}
\newcommand{\IH}{\mathscr{H}}
\newcommand{\hendo}{\mathrm{End}}
\newcommand{\Arr}{\mathcal{A}}
\newcommand{\sym}{\mathfrak{S}}

\newcommand{\young}{\mathscr{Y}}

\newcommand{\card}{\mathrm{card}\,} 

\newcommand{\GL}{\mathrm{GL}}

\newcommand{\GB}{\mathrm{B}}
\newcommand{\For}{\mathbb{F}}
\newcommand{\proba}{\mathbb{P}}
\newcommand{\esper}{\mathbb{E}}
\newcommand{\eps}{\varepsilon}
\newcommand{\lle}{\left[\!\left[} 
\newcommand{\rre}{\right]\!\right]} 
\newcommand{\obs}{\mathscr{O}}
\newcommand{\tilp}{\widetilde{p}}
\newcommand{\id}{\mathrm{id}}
\newcommand{\iso}{\stackrel{\text{\tiny iso}}{\simeq}}
\newcommand{\scal}[2]{\left\langle #1\vphantom{#2}\,\right |\left.#2 \vphantom{#1}\right\rangle}
\newcommand{\figcap}[2]{\begin{figure}[ht] \begin{center} {\footnotesize{#1}} \caption{#2} \end{center} \end{figure}}
\DeclareMathOperator{\imaj}{imaj}
\DeclareMathOperator{\Norm}{Norm}

\author{Valentin F\'eray}
\address{Laboratoire Bordelais de Recherche en Informatique (LaBRI),
Universit\'e Bordeaux 1, 351 cours de la Lib\'eration,
33 400 Talence, France}
\email{feray@labri.fr}

\author{Pierre-Lo\"ic M\'eliot}
\address{The Gaspard--Monge Institut of electronique and computer science,
University of Marne-La-Vall\'ee Paris-Est,
77454 Marne-la-Vall\'ee Cedex 2, France}
\email{meliot@phare.normalesup.org}

\title{Asymptotics of $q$-Plancherel measures}

\begin{document}

\begin{abstract}
 In this paper, we are interested in the asymptotic size of the rows and columns of a random Young diagram under a natural deformation of the Plancherel measure coming from Hecke algebras. The first lines of such diagrams are typically of order $n$, so it does not fit in the context of the work of P. Biane and P. \'Sniady. Using the theory of polynomial functions on Young diagrams of Kerov and Olshanski, we are able to compute explicitly the first- and second-order asymptotics of the length of the first rows. Our method also works for other measures, for example those coming from Schur-Weyl representations.
\end{abstract}

\maketitle

\section{Introduction}
\subsection{Background}\label{subsect:background}
Given a combinatorial family of planar objects, it is natural to ask if, with a good rescaling, there is a limiting continuous object when the size tends to infinity. In the case of Young diagrams endowed with Plancherel measure, this question has been solved in 1977 independently by A. Vershik and S. Kerov \cite{KV77} and B. Logan and L. Shepp \cite{LS77}. V. Ivanov and G. Olshanski \cite{IO02}, using unpublished notes of S. Kerov, have been able to describe the fluctuations of diagrams around this limit shape.\bigskip

The Plancherel measure is natural when we look at Young diagrams of a given size $n$ (the set of which will be denoted $\young_{n}$), because it is related to the regular representation of $\sym_n$. But there are many other interesting measures on Young diagrams and we can look for a limiting continuous object in these cases as well. P. Biane has proven the existence of a limiting shape under some general conditions on the character values of the representation associated to the measure \cite{Bia01}. This includes some natural representations, like the one on $(\C^N)^{\otimes n}$ appearing in Schur-Weyl duality under some conditions on $n$ and $N$. In this context, the second-order asymptotics have been studied by P. \'Sniady \cite{Sniady06}.\bigskip

The results mentioned above give a limit shape for the Young diagrams rescaled by $1/\sqrt{n}$,  where $n$ is the number of boxes of the diagram. In particular, the existence of a limit shape implies that the probability that a random diagram has a long row or column (of length $\gg \sqrt{n}$) tends to $0$. Therefore the same method can not be used directly to study the asymptotic behavior of non-balanced Young diagrams, such as diagrams under the $q$-Plancherel measure.\bigskip

This measure is a natural $q$-analog of the Plancherel measure. It was first introduced by Kerov in \cite{KerovQHookWalk} {\it via} a deformed hook walk algorithm.
The asymptotic behavior of diagrams under this measure has then been considered by a Strahov \cite{Stra08}.
By generalizing Kerov's method of a differential growth model \cite{Kerov96}, he has obtained a limit shape for a diagram under the $q$-Plancherel measure with a parameter $q$ depending of the size of the diagrams.\bigskip

We have mentioned above that, if we fix $q$, there is no hope in finding a limit shape of diagrams rescaled by $1/\sqrt{n}$ as the first rows or columns are not of order  $\Theta(\sqrt{n})$ . In this paper, we give the precise behavior (first- and second-order asymptotics) of the length of these first rows or columns (either rows or columns, depending on the value of $q$, are of order $\Theta(n)$). This is a good description of the shape of the diagram as a finite number of rows or columns contain an arbitrarily big proportion of the boxes of the Young diagram (outside a set of small probability). Our method is close to the one of the papers \cite{IO02,Sniady06}. The first-order asymptotics actually follow from general results of harmonic analysis on the infinite symmetric group (\cite{Tho64,KV81,KOO97,KOV04}); our computations are quite independent of these results, and they allow to tackle with the same formalism the second-order asymptotics.\bigskip

A natural question arising after this result on $q$-Plancherel measure is the following: can our method be adapted to measures for which typical diagrams have rows of intermediate scaling $\Theta(n^\alpha)$ with $1/2 < \alpha < 1$: intermediate between the balanced case ($\alpha=1/2$) and the $q$-Plancherel case ($\alpha=1$)? The answer is yes and, as examples, we look at measures associated to Schur-Weyl duality. In these cases, the method gives results on the global shape of the diagram, but not the precise behavior of the first rows and columns. Indeed, as for balanced diagrams, the asymptotic behavior of the lengths of the first rows and columns is not a direct consequence of the existence of the limit shape. For instance, in the case of the Plancherel measure, a consequent amount of extra work is needed: see \cite{BDJ99firstrow, BDJ99secondrow, Okounkov99, BOO00,JohanssonPlancherel}.\bigskip

\subsection{Results}
\subsubsection{First rows of a random Young diagram under the $q$-Plancherel measure}

Let us suppose $q<1$. In this paper, we prove a law of large numbers and a central limit theorem for the vector containing the lengths of the first $k$ rows of a random diagram under the $q$-Plancherel measure $M_{n,q}$. Note that Theorem \ref{th:first_order_as_rows} implies that the first $k$ rows contain asymptotically $(1-q^k) \times n$ boxes, therefore describing the length of theses rows gives a good description of the diagram. As $M_{n,q^{-1}}(\lambda')=M_{n,q}(\lambda)$ (this follows from Definition \ref{def_measure}), in the case $q>1$, our result gives the asymptotic behavior of the columns, which are, in this case, the relevant quantities.

\begin{theorem}[First order asymptotics]\label{th:first_order_as_rows}
Suppose $q<1$. If $\lambda=(\lambda_1\geq \lambda_2\geq \cdots)$ denotes a Young diagram of size $n$ under the $q$-Plancherel measure (\emph{cf.} figure \ref{fig:ex_diag_qplanch}), then
\begin{align*}
\forall i \geq 1,\,\,\,\frac{\lambda_i}{n} & \longrightarrow_{M_{n,q}} (1-q)\,q^{i-1};
\end{align*}
where the arrow $\to_{M_{n,q}}$ means convergence in probability.
\end{theorem}

\figcap{\psset{unit=1.7mm} \pspicture(0,0)(90,10)\psline[linewidth=0.25pt](0.000000,0.000000)(0.891089,0.000000)(0.891089,0.891089)(0.000000,0.891089)(0.000000,0.000000) \psline[linewidth=0.25pt](0.891089,0.000000)(1.78218,0.000000)(1.78218,0.891089)(0.891089,0.891089)(0.891089,0.000000) \psline[linewidth=0.25pt](1.78218,0.000000)(2.67327,0.000000)(2.67327,0.891089)(1.78218,0.891089)(1.78218,0.000000) \psline[linewidth=0.25pt](2.67327,0.000000)(3.56436,0.000000)(3.56436,0.891089)(2.67327,0.891089)(2.67327,0.000000) \psline[linewidth=0.25pt](3.56436,0.000000)(4.45545,0.000000)(4.45545,0.891089)(3.56436,0.891089)(3.56436,0.000000) \psline[linewidth=0.25pt](4.45545,0.000000)(5.34653,0.000000)(5.34653,0.891089)(4.45545,0.891089)(4.45545,0.000000) \psline[linewidth=0.25pt](5.34653,0.000000)(6.23762,0.000000)(6.23762,0.891089)(5.34653,0.891089)(5.34653,0.000000) \psline[linewidth=0.25pt](6.23762,0.000000)(7.12871,0.000000)(7.12871,0.891089)(6.23762,0.891089)(6.23762,0.000000) \psline[linewidth=0.25pt](7.12871,0.000000)(8.01980,0.000000)(8.01980,0.891089)(7.12871,0.891089)(7.12871,0.000000) \psline[linewidth=0.25pt](8.01980,0.000000)(8.91089,0.000000)(8.91089,0.891089)(8.01980,0.891089)(8.01980,0.000000) \psline[linewidth=0.25pt](8.91089,0.000000)(9.80198,0.000000)(9.80198,0.891089)(8.91089,0.891089)(8.91089,0.000000) \psline[linewidth=0.25pt](9.80198,0.000000)(10.6931,0.000000)(10.6931,0.891089)(9.80198,0.891089)(9.80198,0.000000) \psline[linewidth=0.25pt](10.6931,0.000000)(11.5842,0.000000)(11.5842,0.891089)(10.6931,0.891089)(10.6931,0.000000) \psline[linewidth=0.25pt](11.5842,0.000000)(12.4752,0.000000)(12.4752,0.891089)(11.5842,0.891089)(11.5842,0.000000) \psline[linewidth=0.25pt](12.4752,0.000000)(13.3663,0.000000)(13.3663,0.891089)(12.4752,0.891089)(12.4752,0.000000) \psline[linewidth=0.25pt](13.3663,0.000000)(14.2574,0.000000)(14.2574,0.891089)(13.3663,0.891089)(13.3663,0.000000) \psline[linewidth=0.25pt](14.2574,0.000000)(15.1485,0.000000)(15.1485,0.891089)(14.2574,0.891089)(14.2574,0.000000) \psline[linewidth=0.25pt](15.1485,0.000000)(16.0396,0.000000)(16.0396,0.891089)(15.1485,0.891089)(15.1485,0.000000) \psline[linewidth=0.25pt](16.0396,0.000000)(16.9307,0.000000)(16.9307,0.891089)(16.0396,0.891089)(16.0396,0.000000) \psline[linewidth=0.25pt](16.9307,0.000000)(17.8218,0.000000)(17.8218,0.891089)(16.9307,0.891089)(16.9307,0.000000) \psline[linewidth=0.25pt](17.8218,0.000000)(18.7129,0.000000)(18.7129,0.891089)(17.8218,0.891089)(17.8218,0.000000) \psline[linewidth=0.25pt](18.7129,0.000000)(19.6040,0.000000)(19.6040,0.891089)(18.7129,0.891089)(18.7129,0.000000) \psline[linewidth=0.25pt](19.6040,0.000000)(20.4950,0.000000)(20.4950,0.891089)(19.6040,0.891089)(19.6040,0.000000) \psline[linewidth=0.25pt](20.4950,0.000000)(21.3861,0.000000)(21.3861,0.891089)(20.4950,0.891089)(20.4950,0.000000) \psline[linewidth=0.25pt](21.3861,0.000000)(22.2772,0.000000)(22.2772,0.891089)(21.3861,0.891089)(21.3861,0.000000) \psline[linewidth=0.25pt](22.2772,0.000000)(23.1683,0.000000)(23.1683,0.891089)(22.2772,0.891089)(22.2772,0.000000) \psline[linewidth=0.25pt](23.1683,0.000000)(24.0594,0.000000)(24.0594,0.891089)(23.1683,0.891089)(23.1683,0.000000) \psline[linewidth=0.25pt](24.0594,0.000000)(24.9505,0.000000)(24.9505,0.891089)(24.0594,0.891089)(24.0594,0.000000) \psline[linewidth=0.25pt](24.9505,0.000000)(25.8416,0.000000)(25.8416,0.891089)(24.9505,0.891089)(24.9505,0.000000) \psline[linewidth=0.25pt](25.8416,0.000000)(26.7327,0.000000)(26.7327,0.891089)(25.8416,0.891089)(25.8416,0.000000) \psline[linewidth=0.25pt](26.7327,0.000000)(27.6238,0.000000)(27.6238,0.891089)(26.7327,0.891089)(26.7327,0.000000) \psline[linewidth=0.25pt](27.6238,0.000000)(28.5149,0.000000)(28.5149,0.891089)(27.6238,0.891089)(27.6238,0.000000) \psline[linewidth=0.25pt](28.5149,0.000000)(29.4059,0.000000)(29.4059,0.891089)(28.5149,0.891089)(28.5149,0.000000) \psline[linewidth=0.25pt](29.4059,0.000000)(30.2970,0.000000)(30.2970,0.891089)(29.4059,0.891089)(29.4059,0.000000) \psline[linewidth=0.25pt](30.2970,0.000000)(31.1881,0.000000)(31.1881,0.891089)(30.2970,0.891089)(30.2970,0.000000) \psline[linewidth=0.25pt](31.1881,0.000000)(32.0792,0.000000)(32.0792,0.891089)(31.1881,0.891089)(31.1881,0.000000) \psline[linewidth=0.25pt](32.0792,0.000000)(32.9703,0.000000)(32.9703,0.891089)(32.0792,0.891089)(32.0792,0.000000) \psline[linewidth=0.25pt](32.9703,0.000000)(33.8614,0.000000)(33.8614,0.891089)(32.9703,0.891089)(32.9703,0.000000) \psline[linewidth=0.25pt](33.8614,0.000000)(34.7525,0.000000)(34.7525,0.891089)(33.8614,0.891089)(33.8614,0.000000) \psline[linewidth=0.25pt](34.7525,0.000000)(35.6436,0.000000)(35.6436,0.891089)(34.7525,0.891089)(34.7525,0.000000) \psline[linewidth=0.25pt](35.6436,0.000000)(36.5347,0.000000)(36.5347,0.891089)(35.6436,0.891089)(35.6436,0.000000) \psline[linewidth=0.25pt](36.5347,0.000000)(37.4257,0.000000)(37.4257,0.891089)(36.5347,0.891089)(36.5347,0.000000) \psline[linewidth=0.25pt](37.4257,0.000000)(38.3168,0.000000)(38.3168,0.891089)(37.4257,0.891089)(37.4257,0.000000) \psline[linewidth=0.25pt](38.3168,0.000000)(39.2079,0.000000)(39.2079,0.891089)(38.3168,0.891089)(38.3168,0.000000) \psline[linewidth=0.25pt](39.2079,0.000000)(40.0990,0.000000)(40.0990,0.891089)(39.2079,0.891089)(39.2079,0.000000) \psline[linewidth=0.25pt](40.0990,0.000000)(40.9901,0.000000)(40.9901,0.891089)(40.0990,0.891089)(40.0990,0.000000) \psline[linewidth=0.25pt](40.9901,0.000000)(41.8812,0.000000)(41.8812,0.891089)(40.9901,0.891089)(40.9901,0.000000) \psline[linewidth=0.25pt](41.8812,0.000000)(42.7723,0.000000)(42.7723,0.891089)(41.8812,0.891089)(41.8812,0.000000) \psline[linewidth=0.25pt](42.7723,0.000000)(43.6634,0.000000)(43.6634,0.891089)(42.7723,0.891089)(42.7723,0.000000) \psline[linewidth=0.25pt](43.6634,0.000000)(44.5545,0.000000)(44.5545,0.891089)(43.6634,0.891089)(43.6634,0.000000) \psline[linewidth=0.25pt](44.5545,0.000000)(45.4455,0.000000)(45.4455,0.891089)(44.5545,0.891089)(44.5545,0.000000) \psline[linewidth=0.25pt](45.4455,0.000000)(46.3366,0.000000)(46.3366,0.891089)(45.4455,0.891089)(45.4455,0.000000) \psline[linewidth=0.25pt](46.3366,0.000000)(47.2277,0.000000)(47.2277,0.891089)(46.3366,0.891089)(46.3366,0.000000) \psline[linewidth=0.25pt](47.2277,0.000000)(48.1188,0.000000)(48.1188,0.891089)(47.2277,0.891089)(47.2277,0.000000) \psline[linewidth=0.25pt](48.1188,0.000000)(49.0099,0.000000)(49.0099,0.891089)(48.1188,0.891089)(48.1188,0.000000) \psline[linewidth=0.25pt](49.0099,0.000000)(49.9010,0.000000)(49.9010,0.891089)(49.0099,0.891089)(49.0099,0.000000) \psline[linewidth=0.25pt](49.9010,0.000000)(50.7921,0.000000)(50.7921,0.891089)(49.9010,0.891089)(49.9010,0.000000) \psline[linewidth=0.25pt](50.7921,0.000000)(51.6832,0.000000)(51.6832,0.891089)(50.7921,0.891089)(50.7921,0.000000) \psline[linewidth=0.25pt](51.6832,0.000000)(52.5743,0.000000)(52.5743,0.891089)(51.6832,0.891089)(51.6832,0.000000) \psline[linewidth=0.25pt](52.5743,0.000000)(53.4653,0.000000)(53.4653,0.891089)(52.5743,0.891089)(52.5743,0.000000) \psline[linewidth=0.25pt](53.4653,0.000000)(54.3564,0.000000)(54.3564,0.891089)(53.4653,0.891089)(53.4653,0.000000) \psline[linewidth=0.25pt](54.3564,0.000000)(55.2475,0.000000)(55.2475,0.891089)(54.3564,0.891089)(54.3564,0.000000) \psline[linewidth=0.25pt](55.2475,0.000000)(56.1386,0.000000)(56.1386,0.891089)(55.2475,0.891089)(55.2475,0.000000) \psline[linewidth=0.25pt](56.1386,0.000000)(57.0297,0.000000)(57.0297,0.891089)(56.1386,0.891089)(56.1386,0.000000) \psline[linewidth=0.25pt](57.0297,0.000000)(57.9208,0.000000)(57.9208,0.891089)(57.0297,0.891089)(57.0297,0.000000) \psline[linewidth=0.25pt](57.9208,0.000000)(58.8119,0.000000)(58.8119,0.891089)(57.9208,0.891089)(57.9208,0.000000) \psline[linewidth=0.25pt](58.8119,0.000000)(59.7030,0.000000)(59.7030,0.891089)(58.8119,0.891089)(58.8119,0.000000) \psline[linewidth=0.25pt](59.7030,0.000000)(60.5941,0.000000)(60.5941,0.891089)(59.7030,0.891089)(59.7030,0.000000) \psline[linewidth=0.25pt](60.5941,0.000000)(61.4851,0.000000)(61.4851,0.891089)(60.5941,0.891089)(60.5941,0.000000) \psline[linewidth=0.25pt](61.4851,0.000000)(62.3762,0.000000)(62.3762,0.891089)(61.4851,0.891089)(61.4851,0.000000) \psline[linewidth=0.25pt](62.3762,0.000000)(63.2673,0.000000)(63.2673,0.891089)(62.3762,0.891089)(62.3762,0.000000) \psline[linewidth=0.25pt](63.2673,0.000000)(64.1584,0.000000)(64.1584,0.891089)(63.2673,0.891089)(63.2673,0.000000) \psline[linewidth=0.25pt](64.1584,0.000000)(65.0495,0.000000)(65.0495,0.891089)(64.1584,0.891089)(64.1584,0.000000) \psline[linewidth=0.25pt](65.0495,0.000000)(65.9406,0.000000)(65.9406,0.891089)(65.0495,0.891089)(65.0495,0.000000) \psline[linewidth=0.25pt](65.9406,0.000000)(66.8317,0.000000)(66.8317,0.891089)(65.9406,0.891089)(65.9406,0.000000) \psline[linewidth=0.25pt](66.8317,0.000000)(67.7228,0.000000)(67.7228,0.891089)(66.8317,0.891089)(66.8317,0.000000) \psline[linewidth=0.25pt](67.7228,0.000000)(68.6139,0.000000)(68.6139,0.891089)(67.7228,0.891089)(67.7228,0.000000) \psline[linewidth=0.25pt](68.6139,0.000000)(69.5050,0.000000)(69.5050,0.891089)(68.6139,0.891089)(68.6139,0.000000) \psline[linewidth=0.25pt](69.5050,0.000000)(70.3960,0.000000)(70.3960,0.891089)(69.5050,0.891089)(69.5050,0.000000) \psline[linewidth=0.25pt](70.3960,0.000000)(71.2871,0.000000)(71.2871,0.891089)(70.3960,0.891089)(70.3960,0.000000) \psline[linewidth=0.25pt](71.2871,0.000000)(72.1782,0.000000)(72.1782,0.891089)(71.2871,0.891089)(71.2871,0.000000) \psline[linewidth=0.25pt](72.1782,0.000000)(73.0693,0.000000)(73.0693,0.891089)(72.1782,0.891089)(72.1782,0.000000) \psline[linewidth=0.25pt](73.0693,0.000000)(73.9604,0.000000)(73.9604,0.891089)(73.0693,0.891089)(73.0693,0.000000) \psline[linewidth=0.25pt](73.9604,0.000000)(74.8515,0.000000)(74.8515,0.891089)(73.9604,0.891089)(73.9604,0.000000) \psline[linewidth=0.25pt](74.8515,0.000000)(75.7426,0.000000)(75.7426,0.891089)(74.8515,0.891089)(74.8515,0.000000) \psline[linewidth=0.25pt](75.7426,0.000000)(76.6337,0.000000)(76.6337,0.891089)(75.7426,0.891089)(75.7426,0.000000) \psline[linewidth=0.25pt](76.6337,0.000000)(77.5247,0.000000)(77.5247,0.891089)(76.6337,0.891089)(76.6337,0.000000) \psline[linewidth=0.25pt](77.5247,0.000000)(78.4158,0.000000)(78.4158,0.891089)(77.5247,0.891089)(77.5247,0.000000) \psline[linewidth=0.25pt](78.4158,0.000000)(79.3069,0.000000)(79.3069,0.891089)(78.4158,0.891089)(78.4158,0.000000) \psline[linewidth=0.25pt](79.3069,0.000000)(80.1980,0.000000)(80.1980,0.891089)(79.3069,0.891089)(79.3069,0.000000) \psline[linewidth=0.25pt](80.1980,0.000000)(81.0891,0.000000)(81.0891,0.891089)(80.1980,0.891089)(80.1980,0.000000) \psline[linewidth=0.25pt](81.0891,0.000000)(81.9802,0.000000)(81.9802,0.891089)(81.0891,0.891089)(81.0891,0.000000) \psline[linewidth=0.25pt](81.9802,0.000000)(82.8713,0.000000)(82.8713,0.891089)(81.9802,0.891089)(81.9802,0.000000) \psline[linewidth=0.25pt](82.8713,0.000000)(83.7624,0.000000)(83.7624,0.891089)(82.8713,0.891089)(82.8713,0.000000) \psline[linewidth=0.25pt](83.7624,0.000000)(84.6535,0.000000)(84.6535,0.891089)(83.7624,0.891089)(83.7624,0.000000) \psline[linewidth=0.25pt](84.6535,0.000000)(85.5446,0.000000)(85.5446,0.891089)(84.6535,0.891089)(84.6535,0.000000) \psline[linewidth=0.25pt](85.5446,0.000000)(86.4356,0.000000)(86.4356,0.891089)(85.5446,0.891089)(85.5446,0.000000) \psline[linewidth=0.25pt](86.4356,0.000000)(87.3267,0.000000)(87.3267,0.891089)(86.4356,0.891089)(86.4356,0.000000) \psline[linewidth=0.25pt](87.3267,0.000000)(88.2178,0.000000)(88.2178,0.891089)(87.3267,0.891089)(87.3267,0.000000) \psline[linewidth=0.25pt](88.2178,0.000000)(89.1089,0.000000)(89.1089,0.891089)(88.2178,0.891089)(88.2178,0.000000) \psline[linewidth=0.25pt](89.1089,0.000000)(90.0000,0.000000)(90.0000,0.891089)(89.1089,0.891089)(89.1089,0.000000) \psline[linewidth=0.25pt](0.000000,0.891089)(0.891089,0.891089)(0.891089,1.78218)(0.000000,1.78218)(0.000000,0.891089) \psline[linewidth=0.25pt](0.891089,0.891089)(1.78218,0.891089)(1.78218,1.78218)(0.891089,1.78218)(0.891089,0.891089) \psline[linewidth=0.25pt](1.78218,0.891089)(2.67327,0.891089)(2.67327,1.78218)(1.78218,1.78218)(1.78218,0.891089) \psline[linewidth=0.25pt](2.67327,0.891089)(3.56436,0.891089)(3.56436,1.78218)(2.67327,1.78218)(2.67327,0.891089) \psline[linewidth=0.25pt](3.56436,0.891089)(4.45545,0.891089)(4.45545,1.78218)(3.56436,1.78218)(3.56436,0.891089) \psline[linewidth=0.25pt](4.45545,0.891089)(5.34653,0.891089)(5.34653,1.78218)(4.45545,1.78218)(4.45545,0.891089) \psline[linewidth=0.25pt](5.34653,0.891089)(6.23762,0.891089)(6.23762,1.78218)(5.34653,1.78218)(5.34653,0.891089) \psline[linewidth=0.25pt](6.23762,0.891089)(7.12871,0.891089)(7.12871,1.78218)(6.23762,1.78218)(6.23762,0.891089) \psline[linewidth=0.25pt](7.12871,0.891089)(8.01980,0.891089)(8.01980,1.78218)(7.12871,1.78218)(7.12871,0.891089) \psline[linewidth=0.25pt](8.01980,0.891089)(8.91089,0.891089)(8.91089,1.78218)(8.01980,1.78218)(8.01980,0.891089) \psline[linewidth=0.25pt](8.91089,0.891089)(9.80198,0.891089)(9.80198,1.78218)(8.91089,1.78218)(8.91089,0.891089) \psline[linewidth=0.25pt](9.80198,0.891089)(10.6931,0.891089)(10.6931,1.78218)(9.80198,1.78218)(9.80198,0.891089) \psline[linewidth=0.25pt](10.6931,0.891089)(11.5842,0.891089)(11.5842,1.78218)(10.6931,1.78218)(10.6931,0.891089) \psline[linewidth=0.25pt](11.5842,0.891089)(12.4752,0.891089)(12.4752,1.78218)(11.5842,1.78218)(11.5842,0.891089) \psline[linewidth=0.25pt](12.4752,0.891089)(13.3663,0.891089)(13.3663,1.78218)(12.4752,1.78218)(12.4752,0.891089) \psline[linewidth=0.25pt](13.3663,0.891089)(14.2574,0.891089)(14.2574,1.78218)(13.3663,1.78218)(13.3663,0.891089) \psline[linewidth=0.25pt](14.2574,0.891089)(15.1485,0.891089)(15.1485,1.78218)(14.2574,1.78218)(14.2574,0.891089) \psline[linewidth=0.25pt](15.1485,0.891089)(16.0396,0.891089)(16.0396,1.78218)(15.1485,1.78218)(15.1485,0.891089) \psline[linewidth=0.25pt](16.0396,0.891089)(16.9307,0.891089)(16.9307,1.78218)(16.0396,1.78218)(16.0396,0.891089) \psline[linewidth=0.25pt](16.9307,0.891089)(17.8218,0.891089)(17.8218,1.78218)(16.9307,1.78218)(16.9307,0.891089) \psline[linewidth=0.25pt](17.8218,0.891089)(18.7129,0.891089)(18.7129,1.78218)(17.8218,1.78218)(17.8218,0.891089) \psline[linewidth=0.25pt](18.7129,0.891089)(19.6040,0.891089)(19.6040,1.78218)(18.7129,1.78218)(18.7129,0.891089) \psline[linewidth=0.25pt](19.6040,0.891089)(20.4950,0.891089)(20.4950,1.78218)(19.6040,1.78218)(19.6040,0.891089) \psline[linewidth=0.25pt](20.4950,0.891089)(21.3861,0.891089)(21.3861,1.78218)(20.4950,1.78218)(20.4950,0.891089) \psline[linewidth=0.25pt](21.3861,0.891089)(22.2772,0.891089)(22.2772,1.78218)(21.3861,1.78218)(21.3861,0.891089) \psline[linewidth=0.25pt](22.2772,0.891089)(23.1683,0.891089)(23.1683,1.78218)(22.2772,1.78218)(22.2772,0.891089) \psline[linewidth=0.25pt](23.1683,0.891089)(24.0594,0.891089)(24.0594,1.78218)(23.1683,1.78218)(23.1683,0.891089) \psline[linewidth=0.25pt](24.0594,0.891089)(24.9505,0.891089)(24.9505,1.78218)(24.0594,1.78218)(24.0594,0.891089) \psline[linewidth=0.25pt](24.9505,0.891089)(25.8416,0.891089)(25.8416,1.78218)(24.9505,1.78218)(24.9505,0.891089) \psline[linewidth=0.25pt](25.8416,0.891089)(26.7327,0.891089)(26.7327,1.78218)(25.8416,1.78218)(25.8416,0.891089) \psline[linewidth=0.25pt](26.7327,0.891089)(27.6238,0.891089)(27.6238,1.78218)(26.7327,1.78218)(26.7327,0.891089) \psline[linewidth=0.25pt](27.6238,0.891089)(28.5149,0.891089)(28.5149,1.78218)(27.6238,1.78218)(27.6238,0.891089) \psline[linewidth=0.25pt](28.5149,0.891089)(29.4059,0.891089)(29.4059,1.78218)(28.5149,1.78218)(28.5149,0.891089) \psline[linewidth=0.25pt](29.4059,0.891089)(30.2970,0.891089)(30.2970,1.78218)(29.4059,1.78218)(29.4059,0.891089) \psline[linewidth=0.25pt](30.2970,0.891089)(31.1881,0.891089)(31.1881,1.78218)(30.2970,1.78218)(30.2970,0.891089) \psline[linewidth=0.25pt](31.1881,0.891089)(32.0792,0.891089)(32.0792,1.78218)(31.1881,1.78218)(31.1881,0.891089) \psline[linewidth=0.25pt](32.0792,0.891089)(32.9703,0.891089)(32.9703,1.78218)(32.0792,1.78218)(32.0792,0.891089) \psline[linewidth=0.25pt](32.9703,0.891089)(33.8614,0.891089)(33.8614,1.78218)(32.9703,1.78218)(32.9703,0.891089) \psline[linewidth=0.25pt](33.8614,0.891089)(34.7525,0.891089)(34.7525,1.78218)(33.8614,1.78218)(33.8614,0.891089) \psline[linewidth=0.25pt](34.7525,0.891089)(35.6436,0.891089)(35.6436,1.78218)(34.7525,1.78218)(34.7525,0.891089) \psline[linewidth=0.25pt](35.6436,0.891089)(36.5347,0.891089)(36.5347,1.78218)(35.6436,1.78218)(35.6436,0.891089) \psline[linewidth=0.25pt](36.5347,0.891089)(37.4257,0.891089)(37.4257,1.78218)(36.5347,1.78218)(36.5347,0.891089) \psline[linewidth=0.25pt](37.4257,0.891089)(38.3168,0.891089)(38.3168,1.78218)(37.4257,1.78218)(37.4257,0.891089) \psline[linewidth=0.25pt](38.3168,0.891089)(39.2079,0.891089)(39.2079,1.78218)(38.3168,1.78218)(38.3168,0.891089) \psline[linewidth=0.25pt](39.2079,0.891089)(40.0990,0.891089)(40.0990,1.78218)(39.2079,1.78218)(39.2079,0.891089) \psline[linewidth=0.25pt](40.0990,0.891089)(40.9901,0.891089)(40.9901,1.78218)(40.0990,1.78218)(40.0990,0.891089) \psline[linewidth=0.25pt](40.9901,0.891089)(41.8812,0.891089)(41.8812,1.78218)(40.9901,1.78218)(40.9901,0.891089) \psline[linewidth=0.25pt](41.8812,0.891089)(42.7723,0.891089)(42.7723,1.78218)(41.8812,1.78218)(41.8812,0.891089) \psline[linewidth=0.25pt](42.7723,0.891089)(43.6634,0.891089)(43.6634,1.78218)(42.7723,1.78218)(42.7723,0.891089) \psline[linewidth=0.25pt](43.6634,0.891089)(44.5545,0.891089)(44.5545,1.78218)(43.6634,1.78218)(43.6634,0.891089) \psline[linewidth=0.25pt](44.5545,0.891089)(45.4455,0.891089)(45.4455,1.78218)(44.5545,1.78218)(44.5545,0.891089) \psline[linewidth=0.25pt](0.000000,1.78218)(0.891089,1.78218)(0.891089,2.67327)(0.000000,2.67327)(0.000000,1.78218) \psline[linewidth=0.25pt](0.891089,1.78218)(1.78218,1.78218)(1.78218,2.67327)(0.891089,2.67327)(0.891089,1.78218) \psline[linewidth=0.25pt](1.78218,1.78218)(2.67327,1.78218)(2.67327,2.67327)(1.78218,2.67327)(1.78218,1.78218) \psline[linewidth=0.25pt](2.67327,1.78218)(3.56436,1.78218)(3.56436,2.67327)(2.67327,2.67327)(2.67327,1.78218) \psline[linewidth=0.25pt](3.56436,1.78218)(4.45545,1.78218)(4.45545,2.67327)(3.56436,2.67327)(3.56436,1.78218) \psline[linewidth=0.25pt](4.45545,1.78218)(5.34653,1.78218)(5.34653,2.67327)(4.45545,2.67327)(4.45545,1.78218) \psline[linewidth=0.25pt](5.34653,1.78218)(6.23762,1.78218)(6.23762,2.67327)(5.34653,2.67327)(5.34653,1.78218) \psline[linewidth=0.25pt](6.23762,1.78218)(7.12871,1.78218)(7.12871,2.67327)(6.23762,2.67327)(6.23762,1.78218) \psline[linewidth=0.25pt](7.12871,1.78218)(8.01980,1.78218)(8.01980,2.67327)(7.12871,2.67327)(7.12871,1.78218) \psline[linewidth=0.25pt](8.01980,1.78218)(8.91089,1.78218)(8.91089,2.67327)(8.01980,2.67327)(8.01980,1.78218) \psline[linewidth=0.25pt](8.91089,1.78218)(9.80198,1.78218)(9.80198,2.67327)(8.91089,2.67327)(8.91089,1.78218) \psline[linewidth=0.25pt](9.80198,1.78218)(10.6931,1.78218)(10.6931,2.67327)(9.80198,2.67327)(9.80198,1.78218) \psline[linewidth=0.25pt](10.6931,1.78218)(11.5842,1.78218)(11.5842,2.67327)(10.6931,2.67327)(10.6931,1.78218) \psline[linewidth=0.25pt](11.5842,1.78218)(12.4752,1.78218)(12.4752,2.67327)(11.5842,2.67327)(11.5842,1.78218) \psline[linewidth=0.25pt](12.4752,1.78218)(13.3663,1.78218)(13.3663,2.67327)(12.4752,2.67327)(12.4752,1.78218) \psline[linewidth=0.25pt](13.3663,1.78218)(14.2574,1.78218)(14.2574,2.67327)(13.3663,2.67327)(13.3663,1.78218) \psline[linewidth=0.25pt](14.2574,1.78218)(15.1485,1.78218)(15.1485,2.67327)(14.2574,2.67327)(14.2574,1.78218) \psline[linewidth=0.25pt](15.1485,1.78218)(16.0396,1.78218)(16.0396,2.67327)(15.1485,2.67327)(15.1485,1.78218) \psline[linewidth=0.25pt](16.0396,1.78218)(16.9307,1.78218)(16.9307,2.67327)(16.0396,2.67327)(16.0396,1.78218) \psline[linewidth=0.25pt](16.9307,1.78218)(17.8218,1.78218)(17.8218,2.67327)(16.9307,2.67327)(16.9307,1.78218) \psline[linewidth=0.25pt](17.8218,1.78218)(18.7129,1.78218)(18.7129,2.67327)(17.8218,2.67327)(17.8218,1.78218) \psline[linewidth=0.25pt](18.7129,1.78218)(19.6040,1.78218)(19.6040,2.67327)(18.7129,2.67327)(18.7129,1.78218) \psline[linewidth=0.25pt](19.6040,1.78218)(20.4950,1.78218)(20.4950,2.67327)(19.6040,2.67327)(19.6040,1.78218) \psline[linewidth=0.25pt](20.4950,1.78218)(21.3861,1.78218)(21.3861,2.67327)(20.4950,2.67327)(20.4950,1.78218) \psline[linewidth=0.25pt](21.3861,1.78218)(22.2772,1.78218)(22.2772,2.67327)(21.3861,2.67327)(21.3861,1.78218) \psline[linewidth=0.25pt](22.2772,1.78218)(23.1683,1.78218)(23.1683,2.67327)(22.2772,2.67327)(22.2772,1.78218) \psline[linewidth=0.25pt](23.1683,1.78218)(24.0594,1.78218)(24.0594,2.67327)(23.1683,2.67327)(23.1683,1.78218) \psline[linewidth=0.25pt](24.0594,1.78218)(24.9505,1.78218)(24.9505,2.67327)(24.0594,2.67327)(24.0594,1.78218) \psline[linewidth=0.25pt](0.000000,2.67327)(0.891089,2.67327)(0.891089,3.56436)(0.000000,3.56436)(0.000000,2.67327) \psline[linewidth=0.25pt](0.891089,2.67327)(1.78218,2.67327)(1.78218,3.56436)(0.891089,3.56436)(0.891089,2.67327) \psline[linewidth=0.25pt](1.78218,2.67327)(2.67327,2.67327)(2.67327,3.56436)(1.78218,3.56436)(1.78218,2.67327) \psline[linewidth=0.25pt](2.67327,2.67327)(3.56436,2.67327)(3.56436,3.56436)(2.67327,3.56436)(2.67327,2.67327) \psline[linewidth=0.25pt](3.56436,2.67327)(4.45545,2.67327)(4.45545,3.56436)(3.56436,3.56436)(3.56436,2.67327) \psline[linewidth=0.25pt](4.45545,2.67327)(5.34653,2.67327)(5.34653,3.56436)(4.45545,3.56436)(4.45545,2.67327) \psline[linewidth=0.25pt](5.34653,2.67327)(6.23762,2.67327)(6.23762,3.56436)(5.34653,3.56436)(5.34653,2.67327) \psline[linewidth=0.25pt](6.23762,2.67327)(7.12871,2.67327)(7.12871,3.56436)(6.23762,3.56436)(6.23762,2.67327) \psline[linewidth=0.25pt](0.000000,3.56436)(0.891089,3.56436)(0.891089,4.45545)(0.000000,4.45545)(0.000000,3.56436) \psline[linewidth=0.25pt](0.891089,3.56436)(1.78218,3.56436)(1.78218,4.45545)(0.891089,4.45545)(0.891089,3.56436) \psline[linewidth=0.25pt](1.78218,3.56436)(2.67327,3.56436)(2.67327,4.45545)(1.78218,4.45545)(1.78218,3.56436) \psline[linewidth=0.25pt](2.67327,3.56436)(3.56436,3.56436)(3.56436,4.45545)(2.67327,4.45545)(2.67327,3.56436) \psline[linewidth=0.25pt](3.56436,3.56436)(4.45545,3.56436)(4.45545,4.45545)(3.56436,4.45545)(3.56436,3.56436) \psline[linewidth=0.25pt](4.45545,3.56436)(5.34653,3.56436)(5.34653,4.45545)(4.45545,4.45545)(4.45545,3.56436) \psline[linewidth=0.25pt](5.34653,3.56436)(6.23762,3.56436)(6.23762,4.45545)(5.34653,4.45545)(5.34653,3.56436) \psline[linewidth=0.25pt](0.000000,4.45545)(0.891089,4.45545)(0.891089,5.34653)(0.000000,5.34653)(0.000000,4.45545) \psline[linewidth=0.25pt](0.891089,4.45545)(1.78218,4.45545)(1.78218,5.34653)(0.891089,5.34653)(0.891089,4.45545) \psline[linewidth=0.25pt](1.78218,4.45545)(2.67327,4.45545)(2.67327,5.34653)(1.78218,5.34653)(1.78218,4.45545) \psline[linewidth=0.25pt](0.000000,5.34653)(0.891089,5.34653)(0.891089,6.23762)(0.000000,6.23762)(0.000000,5.34653) \psline[linewidth=0.25pt](0.000000,6.23762)(0.891089,6.23762)(0.891089,7.12871)(0.000000,7.12871)(0.000000,6.23762) \endpspicture}{A 
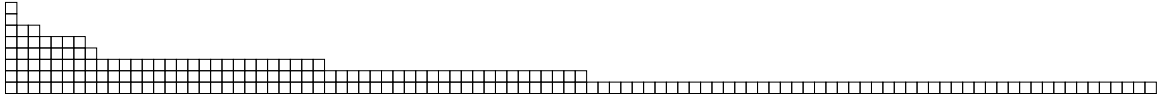
random Young diagram $\lambda=(101,51,28,8,7,3,1,1)$ of size $200$ under the $1/2$-Plancherel measure.\label{fig:ex_diag_qplanch}}

\begin{theorem}[Second order asymptotics]\label{th:second_order_as_rows}
Under the same hypothesis, if $Y_{n,q,i}$ is the rescaled deviation
$$\sqrt{n}\,\left(\frac{\lambda_i}{n} - (1-q)\,q^{i-1}\right),$$
then we have convergence of the finite-dimensional laws of the random process $(Y_{n,q,i})_{i \geq 1}$ towards those of a gaussian process $(Y_{q,i})_{i \geq 1}$ with:
$$\esper[Y_{q,i}]=0 \quad ;\quad \esper[Y_{q,i}^2]=(1-q)\,q^{i-1}-(1-q)^2\,q^{2(i-1)}\quad ;\quad 
\mathrm{cov}(Y_{q,i},Y_{q,j})=-(1-q)^2\,q^{i+j-2}.$$
\end{theorem}\bigskip

\noindent In particular, two different coordinates $Y_{q,i}$ and $Y_{q,j}$ are negatively correlated. For the readers not accustomed to the topology of convergence in law, let us give a more concrete version of Theorem \ref{th:second_order_as_rows}. Given $\eps>0$ and for any positive real numbers $x_1,\ldots,x_r$,  Theorem \ref{th:second_order_as_rows} ensures that, for $n$ big enough:
$$\left|\proba\left[ \forall i \in \lle 1,r \rre,\,\, |\lambda_i- (1-q)\,q^{i-1}\,n|\leq x_i\,\sqrt{n} \right]-\frac{1}{\sqrt{(2\pi)^r\,\det Q}}\,\int_{-x_{1}}^{x_1}\cdots\int_{-x_{r}}^{x_r}\E^{-\frac{{}^tXQ^{-1}X}{2}}\,dX\right| \leq \eps,$$
where $Q$ is the symmetric matrix with coefficients
$Q_{ij}=\delta_{ij}\,(1-q)\,q^{i-1}-(1-q)^2\,q^{i+j-2}$.
This matrix is positive definite because of the Hadamard rule
(in each column, the diagonal coefficient is positive and strictly bigger
than the sum of the absolute values of the other coordinates).\bigskip

\subsubsection{Other measures}

As mentioned in paragraph \ref{subsect:background}, our method is quite general. We can establish limit shape theorems with different rescalings. As an example, in paragraph \ref{sect:schur-weyl}, we introduce some natural probability measures on Young diagrams. We obtain the following result:

\begin{theorem} \label{prop:SW_asymptotic_shape}
Let $(\lambda_n)_{n}$ be the sequence of random Young diagrams described in paragraph \ref{subsect:def_SW}.
We assume $N \sim c\cdot n^\alpha$ with $\alpha <1/2$.
%If we rescale $\lambda_n$ by $n^{1-\alpha}$ in the horizontal direction and $n^{\alpha}$ in the vertical direction,
%the new sequence converges in probability towards a rectangle of side $c$ and $c^{-1}$.\bigskip

Then, for each $\varepsilon, \eta > 0$, there exists an integer $n_0$ such that, for each $n \geq n_0$,
the border of the diagram $\lambda_n$ is, after rescaling and with probability greater than $1 - \varepsilon$, contained
in the hatched area: 
\figcap{\psset{unit=1mm}\pspicture(-2,-2)(100,45)
\pspolygon[fillstyle=hlines](0,25)(0,22)(58,22)(58,0)(90,0)(90,3)(62,3)(62,25)
\psline[linecolor=white,linewidth=2pt](90,0)(90,3)
\psline{->}(-2,0)(95,0)
\psline{->}(0,-2)(0,40)
\rput(105,0){$\times n^{1-\alpha}$}
\rput(0,45){$\times n^{\alpha}$}
\rput(60,-3){$c$}
\rput(-4,25){$c^{-1}$}
\psline{[-]}(30,25)(30,22)
\rput(30,28){$c^{-1}\eta$}
\psline{[-]}(58,13)(62,13)
\rput(66,13){$c\eta$}
\psline(-1,25)(0,25)
\psline(60,-1)(60,1)
\endpspicture
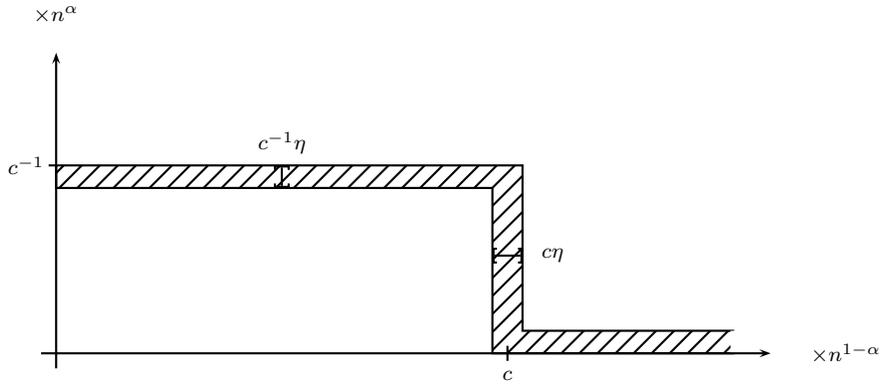
}{Limit shape of a rescaled random diagram under the Schur-Weyl measure of parameters $c>0$ and $\alpha<1/2$.\label{fig:rectangle}}
\end{theorem}

Unfortunately, one can not use this to give the precise behavior of the first row of the diagram. The only result we have is the following (which is a consequence of Proposition \ref{prop:cv_powersums_encore}):

\begin{proposition}\label{prop:SW_first_row}
Let $\lambda_n$ be the sequence of random diagrams described in paragraph \ref{subsect:def_SW}. Let us assume that $N \sim c \, n^\alpha$ with $\alpha <1/2$. For any $\varepsilon > 0$, one has $\lambda_1=o(n^{1-\alpha+\varepsilon})$ in probability.
\end{proposition}
\bigskip

\subsubsection{Longest increasing subsequence} \label{subsect:ss_suite_croissante}

It is well-known that the length of the first row of a Young diagram of size $n$ under Plancherel measure has the same distribution as the length of the longest increasing subsequence of a random permutation of size $n$ (under the uniform distribution). Indeed, the shape of the tableaux obtained by the Robinson-Schensted bijection defines a map
$$\sym_n \longrightarrow \{\lambda \vdash n\},$$
which sends the length of the longest increasing subsequence on the length of the first row of $\lambda$.
The image of the uniform distribution on $\sym_n$ is the Plancherel measure on Young diagrams of size $n$, and this explains the aforementioned result.\bigskip

Our results on the $q$-Plancherel measure can also be interpreted in the context of longest increasing subsequences. Indeed, E. Strahov \cite{Stra08} has noticed that the $q$-Plancherel measure is the image by the map above of the probability measure defined by
$$\forall \ \sigma,\ \proba[\sigma ] = Z \cdot q^{\imaj(\sigma)},$$
where $\imaj(\sigma)$, called major index, is the sum of the integers $i$ between $1$ and $n-1$ such that $\sigma(i)>\sigma(i+1)$ and $Z=1/\{n!\}_q$ is the normalization constant such that the probabilities sum to $1$. So Theorems \ref{th:first_order_as_rows} and \ref{th:second_order_as_rows} imply the following:

\begin{theorem}
 Let $q$ be a real in $]0,1[$ and $n$ an integer. Consider the distribution on permutations such that the probability of picking a permutation $\sigma$ is proportional to $q^{\imaj(\sigma)}$. The length $\ell(\sigma)$ of the longest increasing subsequence of a permutation under this distribution is equivalent in probability to $n(1-q)$. Moreover, the rescaled deviation $(\ell(\sigma) - n(1-q))/\sqrt{n}$ converges in law towards a centered normal law of variance $q(1-q)$.
\end{theorem}\medskip

As mentioned in subsection \ref{subsect:def_SW}, the measure considered in section \ref{sect:schur-weyl} is the image by the RSK algorithm of the uniform distribution on words of length $n$ with letters from $1$ to $N$. So, Proposition \ref{prop:SW_first_row} implies:
\begin{corollary}\label{th:SW_lis_upbound}
 Fix $c>0$ and $0<\alpha <1/2$. Let $\ell$ be the length of the longest non-decreasing subsequence of a random word with $n$ letters between $1$ and $c \cdot \lfloor n^\alpha \rfloor$. When $n$ tends to infinity, one has in probability: $\ell = o(n^{1-\alpha+\varepsilon})$.
\end{corollary}\bigskip

\subsection{Tools}
Since we are looking at measures on diagrams of a given size, irreducible character values on a given permutation can be seen as random variables. Moreover, in our context, their expectations are easy to compute. Therefore we would like to express other random variables, linked directly to the shape of the diagrams in terms of the character value. This can be done thanks to Kerov's and Olshanski's theory of polynomial functions on the set of Young diagrams (see \cite{KO94}). Note that the guideline of the proof is roughly the same as in Kerov's proof of the central limit theorem for Young diagram under Plancherel's measure (see \cite{IO02}).\bigskip

In this article, we shall use the power sums of Frobenius coordinates. With an adequate gradation on the functions on Young diagrams, we shall easily find their asymptotic behavior under the $q$-Plancherel measure. Going back to the behavior of the first rows is then a technical step.\bigskip

\subsection{Some open questions}
\subsubsection{Precise behavior of the first line}
Proposition \ref{prop:SW_first_row} (respectively Corollary \ref{th:SW_lis_upbound}) gives only an upper bound for the behavior of the first row of the diagrams (respectively the length of the longest non-decreasing subsequence of a random word). It would be nice to obtain the precise behavior (equivalent and perhaps fluctuation) of these random variables. If we follow \cite{Okounkov99}, the first step is to study the asymptotic behavior of the observable $p_k$ with a parameter $k$ which also tends to infinity (as quickly as some power of $n$ for example). Such a result could be deduced from the expression of $p_\rho$ in terms of $\varSigma$'s (normalized characters), but unfortunately, equation (\ref{eq:p_fct_sigma}) is hard to invert.

\subsubsection{Concrete realization of power sums in the center of the symmetric group algebras}
We will see in section \ref{sect:observables} that the relations between some functions on Young diagrams can be understood as relations in the center of the symmetric group algebras via the abstract Fourier transform. Therefore, if we want to understand better the expression of the $p_\rho$'s in terms of the $\varSigma$'s (and we would like to, see the first open problem!), one way to do it would be to describe in a nice way the image of the $p_\rho$'s by the abstract Fourier transform. One could also use the power sums of Jucys-Murphy elements which have a simple definition in $\C\sym_n$ and have the same top  homogeneous component as the $p_\rho$ (with a shift of index and up to a multiplicative constant).

\subsubsection{Approximation by a sum of independant random vectors}
It is interesting to note that the covariance matrix appearing in Theorem \ref{th:second_order_as_rows} is exactly the same as the one obtained by A. Borodin in a very different context (the limiting distribution of the sizes of the blocks of the Jordan form of a random upper triangular nilpotent matrix over a finite field, see \cite{Borodin:limitJordan}). In his paper, A. Borodin found this covariance matrix by approximating his random partition by a sum of independent random vectors. This suggests that the $q$-Plancherel measure could perhaps be approximated the same way. Such a result could be obtained by using the limit form of the diagrams (Theorem \ref{th:first_order_as_rows}) and computing the error on the corresponding transition probabilities. Unfortunately, the transition probabilities for the $q$-Plancherel measure (see Figure \ref{qplancherelprocess}) are hard to deal with and we do not know if this method would succeed in our context.
\bigskip

\subsection{Outline of the paper}
In section \ref{sect:def_q-planch}, we define the $q$-Plancherel measure and recall its basic pro\-perties. Then, in section \ref{sect:observables}, we introduce families of observables of diagrams, which will be useful in the proof of the main theorems. Sections \ref{sect:first_order} and \ref{sect:second_order} are respectively devoted to the proofs of Theorems \ref{th:first_order_as_rows} and \ref{th:second_order_as_rows}.
In section \ref{sect:schur-weyl}, we study other measures with the same tools.
\bigskip

\section{The $q$-Plancherel measure}\label{sect:def_q-planch}

In this section, we define and present basic properties of the $q$-Plancherel measure. The characterization given in paragraph \ref{subsect:qPlancherel_trace} is the one which will be used in the next sections. The algebraic interpretation of paragraph \ref{subsect:algebraic_int_regular_trace} gives some motivation to consider at this measure.\bigskip

Let us fix some notations related to the representation theory of $\sym_{n}$. If $\lambda$ is a partition in $\young_{n}$, it corresponds to an irreducible $\sym_{n}$-module $S^{\lambda}$ --- the Specht module of type $\lambda$ --- whose dimension will be denoted $\dim \lambda$. We shall also denote $\varsigma^{\lambda}$ the irreducible character of $S^{\lambda}$, and $\chi^{\lambda}$ the \emph{normalized} character, so that $\varsigma^{\lambda}=(\dim \lambda)\,\chi^{\lambda}$.\bigskip

\subsection{Definition}

Let us recall the definition of Young graph (\emph{cf.} for example \cite[p. 5]{KOO97}): its set of vertices is the disjoint union $\young=\bigsqcup_{n \in \N} \young_{n}$, and its edges are the pairs $(\lambda,\Lambda)$, where the Young diagram $\Lambda$ can be obtained by adding one box to the Young diagram $\lambda$ (this situation will be denote $\lambda \nearrow \Lambda$). The (standard) Plancherel growth process is the markovian partition-valued process $(\lambda_{n})_{n \in \N}$ defined by:
$$\lambda_{0}=\emptyset\qquad;\qquad p(\lambda,\Lambda)=\begin{cases} \frac{1}{|\Lambda|}\,\frac{\dim \Lambda}{\dim \lambda}& \text{if }\lambda \nearrow \Lambda,\\
0&\text{otherwise}.\end{cases}$$ 
Because of the branching rules for the symmetric groups, $p(\cdot,\cdot)$ is indeed a transition function on $\young$, and the law of the $n$-th partition $\lambda_{n}$ is the Plancherel measure on $\young_{n}$:
$$\proba[\lambda_{n}=\lambda]=\mathbb{1}_{|\lambda|=n}\,\frac{(\dim \lambda)^{2}}{n!}.$$
The $q$-Plancherel measures and the $q$-Plancherel process are natural $q$-analogs
of these objects that have been studied by Kerov \cite{KerovQHookWalk} and more recently
by Strahov \cite{Stra08}. Recall that the dimension $\dim \lambda$ is given by
the hook length formula ($\dim \lambda = n!/\prod_{\oblong \in \lambda} h(\oblong)$,
the hook length $h(\oblong)$ of a box $\oblong=(i,j)$ being $(\lambda_{j}-i)+(\lambda_{i}'-j)+1$).
As a consequence, assuming that $\lambda \nearrow \Lambda$, the transition probability $p(\lambda,\Lambda)$ is a quotient of hook products: $$p(\lambda,\Lambda)=\frac{\prod_{\oblong \in \lambda}h(\oblong) }{ \prod_{\oblong \in \Lambda} h(\oblong)}.$$
S. Kerov has shown that one can deform the markovian process by considering $q$-transition probabilities:
$$p_{q}(\lambda,\Lambda) = q^{b(\Lambda)-b(\lambda)}\,\frac{\prod_{\oblong \in \lambda} \{h(\oblong)\}_{q}}{\prod_{\oblong \in \Lambda} \{h(\oblong)\}_{q}},$$
where $b(\lambda)=\sum_{i=1}^{\ell(\lambda)} (i-1)\lambda_i$ and the $q$-analog $\{n\}_{q}$
of a positive integer is by definition $\frac{q^{n}-1}{q-1}$.
\figcap{\footnotesize{
\psset{unit=1mm}\pspicture(-50,-105)(70,5)
\psline(-2.5,0)(2.5,0)(2.5,5)(-2.5,5)(-2.5,0)
\rput(16,2){$ D=1\,\,;\,\,M_{1,q}=1$}
\psline[linecolor=red]{->}(0,-1.5)(20,-13.5)
\rput(16,-7){\textcolor{red}{$\frac{1}{q+1}$}}
\psline[linecolor=red]{->}(0,-1.5)(-20,-11.5)
\rput(-18,-7){\textcolor{red}{$\frac{q}{q+1}$}}
\psline(15,-20)(25,-20)(25,-15)(15,-15)(15,-20) \psline(20,-20)(20,-15)
\rput(44,-18){$D=1 \,\,;\,\,M_{2,q}=1/(q+1) $}
\psline(-17.5,-12.5)(-22.5,-12.5)(-22.5,-22.5)(-17.5,-22.5)(-17.5,-12.5) \psline(-17.5,-17.5)(-22.5,-17.5)
\rput(-35.5,-19){$M_{2,q}=q/(q+1)$}
\rput(-39,-15){$D=q$}
\psline(-5,-35)(5,-35)(5,-40)(0,-40)(0,-45)(-5,-45)(-5,-35) \psline(0,-35)(0,-40)(-5,-40)
\psline[linecolor=red]{->}(-20,-24)(-40,-31)
\rput(-42,-27){\textcolor{red}{$\frac{q^2}{q^2+q+1}$}}
\psline[linecolor=red]{->}(-20,-24)(-1,-34)
\rput(-8,-25){\textcolor{red}{$\frac{q+1}{q^2+q+1}$}}
\psline[linecolor=red]{->}(20,-21.5)(1,-34)
\rput(7,-24.5){\textcolor{red}{$\frac{q^2+q}{q^2+q+1}$}}
\psline[linecolor=red]{->}(20,-21.5)(37.5,-34)
\rput(37,-28){\textcolor{red}{$\frac{1}{q^2+q+1}$}}
\rput(16,-47){$M_{3,q}=(2q)/(q^2+q+1)$}
\rput(10.5,-43){$D=q^2+q$}
\psline[linecolor=red]{->}(-2.5,-46.5)(-2.5,-58.5)
\rput(10,-53){\textcolor{red}{$\frac{q^3+q}{q^4+2q^3+2q^2+2q+1}$}}
\psline(-37.5,-32.5)(-37.5,-47.5)(-42.5,-47.5)(-42.5,-32.5)(-37.5,-32.5) \psline(-37.5,-37.5)(-42.5,-37.5) \psline(-37.5,-42.5)(-42.5,-42.5)
\rput(-32,-46){$D=q^3$}
\rput(-26,-50){$M_{3,q}=q^3/(q^3+2q^2+2q+1)$}
\psline[linecolor=red]{->}(-43.5,-40)(-50,-59)
\rput(-55.5,-50){\textcolor{red}{$\frac{q^3}{q^3+q^2+q+1}$}}
\psline(-47.5,-60)(-52.5,-60)(-52.5,-80)(-47.5,-80)(-47.5,-60)
\psline(-47.5,-65)(-52.5,-65) \psline(-47.5,-70)(-52.5,-70) \psline(-47.5,-75)(-52.5,-75)
\rput(-58.5,-79){$D=q^6$}
\rput(-50,-83){$M_{4,q}=q^6/(q^6+3q^5+5q^4$}
\rput(-47,-86.5){$+6q^3+5q^2+3q+1)$}
\psline(30,-35)(45,-35)(45,-40)(30,-40)(30,-35) \psline(35,-35)(35,-40) \psline(40,-35)(40,-40)
\rput(65,-40){$M_{3,q}=1/(q^3+2q^2+2q+1)$}
\rput(54,-36){$D=1$}
\psline(-5,-60)(5,-60)(5,-70)(-5,-70)(-5,-60) \psline(-5,-65)(5,-65) \psline(0,-60)(0,-70)
\rput(14,-69){$D=q^4+q^2$}
\psline(-15,-80)(-25,-80)(-25,-95)(-20,-95)(-20,-85)(-15,-85)(-15,-80) \psline(-20,-80)(-20,-85)(-25,-85) \psline(-20,-90)(-25,-90)
\psline[linecolor=red]{->}(-2.5,-46.5)(-19,-79)
\psline[linecolor=red]{->}(-2.5,-46.5)(-10,-75)(22,-84)
\psline[linecolor=red]{->}(-45,-52)(-45.5,-56)(-21,-79)
\psframe*[linecolor=white,fillcolor=white](-40,-75)(-20,-69)
\rput(-30,-72){\textcolor{red}{$\frac{q^2+q+1}{q^3+q^2+q+1}$}}
\psline[linecolor=red](-43.5,-40)(-44.5,-48)
\psframe*[linecolor=white,fillcolor=white](-28,-63.5)(-8.5,-57)
\rput(-18,-60){\textcolor{red}{$\frac{q^4+q^3+q^2}{q^4+2q^3+2q^2+2q+1}$}}
\psline(15,-85)(30,-85)(30,-90)(20,-90)(20,-95)(15,-95)(15,-85) \psline(20,-85)(20,-90)(15,-90) \psline(25,-85)(25,-90)
\rput(32,-93){$D=q^3+q^2+q$}
\rput(43,-97){$M_{4,q}=(3q)/(q^4+2q^3+2q^2+2q+1)$}
\rput(-24,-98){$D=q^5+q^4+q^3$}
\rput(-13.5,-102){$M_{4,q}=(3q^3)/(q^4+2q^3+2q^2+2q+1)$}
\psline(35,-60)(55,-60)(55,-65)(35,-65)(35,-60) \psline(40,-60)(40,-65) \psline(45,-60)(45,-65) \psline(50,-60)(50,-65)
\psline[linecolor=red]{->}(37.5,-41)(45,-59)
\psline[linecolor=red]{->}(37.5,-41)(23.5,-84)
\rput(50,-50){\textcolor{red}{$\frac{1}{q^3+q^2+q+1}$}}
\rput(63,-61){$D=1$}
\rput(72,-65){$M_{4,q}=1/(q^6+3q^5+5q^4$}
\rput(74,-68.5){$+6q^3+5q^2+3q+1)$}
\psframe*[linecolor=white,fillcolor=white](25,-75)(35,-71)
\rput(27,-73){$M_{4,q}=(2q^2)/(q^4+3q^3+4q^2+3q+1)$}
\psframe*[linecolor=white,fillcolor=white](25,-60)(35,-54)
\rput(29,-57){\textcolor{red}{$\frac{q^3+q^2+q}{q^3+q^2+q+1}$}}
\rput(0,-82){\textcolor{red}{$\frac{q^2+q+1}{q^4+2q^3+2q^2+2q+1}$}}
\endpspicture}
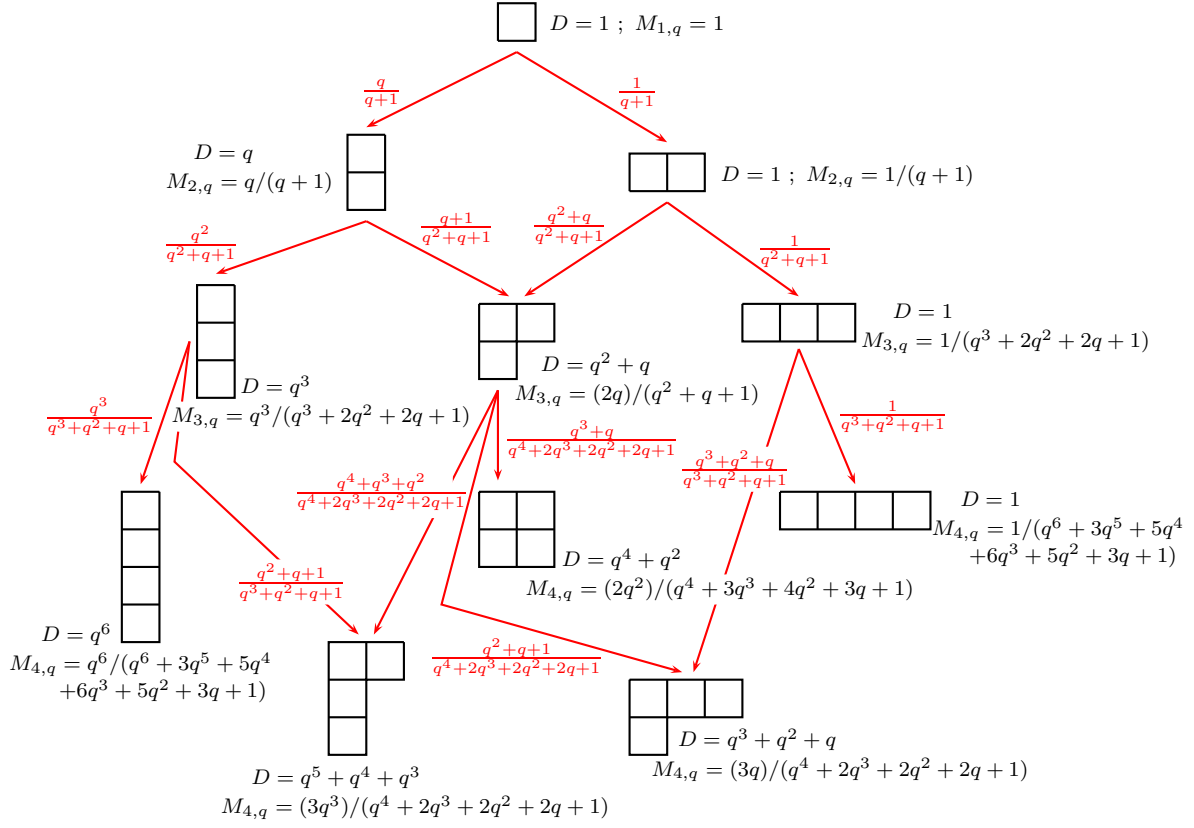
}{Generic degrees, $q$-Plancherel measures and $q$-transition probabilities on diagrams of size $n\leq 4$.\label{qplancherelprocess}}
\begin{definitionproposition}\label{def_measure}
The $q$-Plancherel process is the markovian partition-valued process
$(\lambda_{n})_{n \in \N}$ defined by:
$$\lambda_{0}=\emptyset\qquad;\qquad p_{q}(\lambda,\Lambda)=\begin{cases} q^{b(\Lambda)-b(\lambda)}\,\frac{\prod_{\oblong \in \lambda} \{h(\oblong)\}_{q}}{\prod_{\oblong \in \Lambda} \{h(\oblong)\}_{q}}& \text{if }\lambda \nearrow \Lambda \text{ in the Young graph},\\
0&\text{otherwise}.\end{cases}$$ 
Its marginal distribution of the $n$-th partition of the $q$-Plancherel process is called the $q$-Plancherel measure.
This probability measure $M_{n,q}$ on the set $\young_n$ of Young diagrams of size $n$ can be alternatively defined by:
$$\forall \lambda \vdash n,\,\,\, M_{n,q}(\{\lambda\})=  M_{n,q}(\lambda)=\frac{D_\lambda(q)\,\dim \lambda}{\{n!\}_q},$$
where $\dim \lambda$ is the dimension of the irreducible representation of $\sym_n$ indexed by $\lambda$, and $D_\lambda(q)$ is its generic degree as defined in \cite[Chapter 8]{GP00}:
$$D_{\lambda}(q)=q^{b(\lambda)}\,\frac{\{n!\}_{q}}{\prod_{\oblong \in \lambda}\{h(\oblong)\}_{q}}.$$
\end{definitionproposition}
\noindent The equivalence between the two definitions (as marginal distribution of a markovian process or {\it via} a direct formula for the probability to pick each diagram) is explained in \cite[Section 3]{Stra08}.

In the next two paragraphs, we explain why this $q$-analog is natural in the context of Hecke algebras and general linear groups in finite characteristic.\bigskip

\subsection{Link with representation of the Hecke algebras}\label{subsect:qPlancherel_trace}
The $q$-Plancherel measure appears naturally in the abstract harmonic analysis of the Iwahori-Hecke algebras $\IH_{n,q}$. Recall that if $n \in \N$ and $q$ is a complex parameter, then $\IH_{n,q}$ is the complex associative algebra with generators $T_1,T_2,\ldots,T_{n-1}$ and relations:
\begin{align*}
 \text{braid relations:}&\quad  \forall i,\,\,\,T_iT_{i+1}T_i=T_{i+1}T_iT_{i+1} ;\\
\text{commutation relations:}&\quad \forall\, |j-i|\geq 2,\,\,\,T_iT_j=T_jT_i ;\\
\text{quadratic relations:}&\quad  \forall i,\,\,\,(T_i-q)(T_i+1)=0.
\end{align*}
We recover the symmetric group algebra $\C\sym_n$ for $q=1$. For any value of $q$, one can define an element corresponding to a permutation $\sigma \in \sym_n$:
$$T_\sigma=T_{i_1}T_{i_2}\cdots T_{i_k}, $$
where $\sigma=s_{i_1}s_{i_2}\cdots s_{i_k}$ is a minimal decomposition of $\sigma$ in elementary transpositions (thanks to Matsumoto's theorem, the result does not depend on which minimal decomposition we choose). It can be shown that $(T_\sigma)_{\sigma \in \sym_n}$ forms a linear basis of $\IH_{n,q}$ (see for example \cite{Mat99} or \cite{GP00}). Moreover, if $q \neq 0$ is not a root of unity (in the following, we assume that we are in this case, as $q$ is a positive real number), then $\IH_{n,q}$ is a semi-simple algebra abstractly isomorphic to $\C\sym_n$, and consequently has the same theory of representations: the irreducible $\IH_{n,q}$-modules are parameterized by the partitions of size $n$, and each $q$-Specht module $S^\lambda(q)$ has dimension $\dim S^\lambda(q)=\dim \lambda$. We will denote by $\varsigma^\lambda(q)$ the trace of the simple module $S^\lambda(q)$, and by $\chi^{\lambda}(q)$  the \emph{normalized} character, so that $\varsigma^\lambda(q)=(\dim \lambda)\,\chi^\lambda(q)$. Both objects are central linear forms on $\IH_{n,q}$.\bigskip

Let us consider the regular trace of the Iwahori-Hecke algebra:
$$\tau_q(T_\sigma)=\begin{cases}
                  1&\text{if }\sigma=\mathrm{id},\\
                 0&\text{otherwise}.
                 \end{cases}
$$
The decomposition of $\tau_q$ in the basis of normalized $q$-characters $\chi^\lambda(q)$ writes as
(this corresponds to \cite[equation (16)]{Stra08}):
\begin{align}
\tau_q &= \sum_{\lambda \in \young_n} M_{n,q}(\lambda)\,\chi^\lambda(q)=\sum_{\lambda\in \young_n} \frac{D_\lambda(q)}{\{n!\}_q}\,\varsigma^\lambda(q). \label{eq:tau_qPlancherel}
\end{align}
Of course, this characterizes the $q$-Plancherel measure. Therefore, if we consider the $q$-characters on a fixed element $T_\sigma$ as a function of the indexing partition $\lambda$, its expectation  under the $q$-Plancherel measure is $\tau_q(T_\sigma)$, that is either $0$ or $1$. Using this, we will be able to compute expectations of observables of diagrams. 

\begin{remark} We will prove in this article that we have in fact a concentration of the $q$-characters around the regular trace and compute explicitly their fluctuations, see Corollary \ref{corol:concentration_q-characters} and Remark \ref{RemFluctuationsQCharacters}.\bigskip
\end{remark}

\subsection{Interpretation in the representation theory of $\GL_n(\For_q)$}\label{subsect:algebraic_int_regular_trace}
If $q=p^e$ is a power of a prime number, the previous identities may be interpreted in the setting of the representation theory of the group of invertible $n \times n$ matrices $\GL(n,\For_q)$. Indeed, it is a well-known fact \cite{Iwa64} that the relations defining $\IH_{n,q}$ are those of the convolution algebra
$$\C[\GB(n,\For_q) \backslash \GL(n,\For_q)/ \GB(n,\For_q)] = \left\{f \in \C[\GL(n,\For_q)]\,\,|\,\,\forall b,b' \in \GB(n,\For_q),\,\, f(bgb')=f(g) \right\},$$
where $\GB(n,\For_q)$ is the Borel subgroup of $\GL(n,\For_q)$, containing the invertible upper-triangular matrices.\bigskip

Now, consider the $\GL(n,\For_q)$-module of functions on the flag variety $\GL(n,\For_q)/ \GB(n,\For_q)$ --- or, in other words, right $\GB(n,\For_q)$-invariant complex functions on $\GL(n,\For_q)$. Since $\GL(n,\For_q)$ acts transitively on the flag variety, a $\GL(n,\For_q)$-endomorphism $u$ of this module is entirely determined by a function $h \in \C[\GL(n,\For_q)]$ such that:
$$u(\mathbb{1}_{\GB(n,\For_q)})=\sum_{g \in G} h(g)\,\mathbb{1}_g.$$
Moreover, the condition $u \in \hendo_{\GL(n,\For_q)}\, \C[\GL(n,\For_q)/\GB(n,\For_q)]$ implies that $h$ is left and right $\GB(n,\For_q)$-invariant. As a consequence, we have an isomorphism
$$\IH_{n,q} \iso \C[\GB(n,\For_q)\backslash \GL(n,\For_q) / \GB(n,\For_q)] \iso \hendo_{\GL(n,\For_q)}\,\C[\GL(n,\For_q)/\GB(n,\For_q)],$$
whence a decomposition of the $(\GL(n,\For_q),\IH_{n,q})$-bimodule $\C[\GL(n,\For_q)/\GB(n,\For_q)]$ similar to the classical Schur-Weyl duality between $\GL(n,\C)$ and $\sym_n$:
$$_{\GL(n,\For_q) \curvearrowright }\big\{\C[\GL(n,\For_q)/\GB(n,\For_q)]\big\}_{\curvearrowleft \IH_{n,q}} =\sum_{\lambda \in \young_n} U_\lambda(q) \otimes_{\C} S^\lambda(q).$$
The $U_\lambda(q)$ are called unipotent modules of $\GL(n,\For_q)$. Each $U_\lambda(q)$ is irreducible and of complex dimension
$$\dim U_\lambda(q) = D_\lambda(q) =q^{b(\lambda)} \, \frac{\{n!\}_q}{\prod_{\oblong \in \lambda} \{h(\oblong)\}_q}.$$
Thus, the natural probability measure $M_{n,q}$ on the set of characters of $\IH_{n,q}$ is also the natural probability measure on the set of unipotent characters of $\GL(n,\For_q)$. Moreover, the regular trace $\tau_{q}$ is just the restriction of the normalized trace of $\hendo_{\C}\, \C[\GL(n,\For_{q})/\GB(n,\For_{q})]$ to the subalgebra of $\GL(n,\For_{q})$-morphisms, that is to say, $\IH_{n,q}$.\bigskip

\begin{remark}
In a much more general context, given a (connected, non-twisted) Chevalley group $G^{F}$ over $\For_{q}$ and a cuspidal irreducible character $\rho$ of a Levi subgroup $L^{F}$, the decomposition in irreducibles of the $G^{F}$-module $\mathrm{R}_{L^{F}}^{G^{F}}(\rho)$ obtained by parabolic induction yields a probability measure on the set of characters of a multiple-parameter Hecke algebra of the ramification group
$$W(\rho)=\{w \in \Norm(L^{F})/L^{F}\,\,|\,\,\rho \circ w = \rho\}.$$
This ramification group is a Coxeter group and a subgroup of the relative Weyl group $N(L^{F})/L^{F}$; moreover, the multiple-parameter Hecke algebra is in fact isomorphic to the group algebra $\C W(\rho)$ (see \cite{HL80,Lus84,Geck93}), so the probability measure is in fact over the isomorphism classes of irreducible representations of $W(\rho)$. The case we are studying in this paper is when $G^{F}=\GL(n,\For_{q})$, $L^{F}=(\For_{q})^{n}$, $\rho=1$ is the trivial character and $W(\rho)=W=\sym_{n}$.
\end{remark}
\bigskip

\section{Observables of diagrams}\label{sect:observables}

In this section, we introduce the algebra of observables of diagrams $\obs$. This subalgebra of the algebra of functions on the set of all Young diagrams $\young=\bigsqcup_{n \in \N} \young_n$ is isomorphic to the algebra of (complex) symmetric functions, and it is known to be a powerful tool in the setting of asymptotic representation theory of symmetric groups, see \cite{IO02}. In this paper, the authors introduced four algebraic bases of $\obs$: the power sums $(p_k)_{k\geq 1}$ in the Frobenius coordinates, the power sums $(\tilp_k)_{k \geq 2}$ in the interlacing coordinates, the normalized characters $(\varSigma_k)_{k \geq 1}$ and the free cumulants $(R_k)_{k \geq 2}$. In the following, we shall use the power sums $p_k$, the normalized characters $\varSigma_k$ and their quantizations $\varSigma_{k,q}$; this section is devoted to their presentation.\bigskip

\subsection{Power sums of Frobenius coordinates}
If $\lambda$ is a partition of size $n$, we recall that its Frobenius coordinates $(a_i,b_i)_{1 \leq i \leq d}$ are defined by $a_i=\lambda_i-i$ and $b_i={\lambda_i}'-i$, with $i$ less than the size $d$ of the diagonal of the diagram. The modified Frobenius coordinates are $(A_\lambda,\overline{B}_\lambda)=\big((a_i^*)_{1 \leq i \leq d},(-b_i^*)_{1 \leq i \leq d}\big)$ with $a_i^*=a_i+1/2$ and $b_i^*=b_i+1/2$, see figure \ref{killinginthename} on the next page. Let us define the function $p_k$ on the set $\young$ of all Young diagrams: $p_k(\lambda)$ is the $k$-th power sum of the alphabet $A_\lambda-\overline{B}_\lambda$ (with the notations of $\lambda$-rings). In other words,
$$p_k(\lambda)=\sum_{i=1}^d \,(a_i^*)^k - (-b_i^*)^k.$$
Note that $p_1(\lambda)$ is simply the size $|\lambda|$ of the diagram (whence the modification of the coordinates).

\figcap{\psset{unit=1mm} \pspicture(-50,-3)(50,50)
\psline[linewidth=0.25pt](40,40)(0,0)(-40,40)(-37.5,42.5)(-30,35)(-27.5,37.5)(-25,35)(-22.5,37.5)(-12.5,27.5)(-7.5,32.5)(-5,30)(-2.5,32.5)(2.5,27.5)(5,30)(10,25)(20,35)(22.5,32.5)(25,35)(27.5,32.5)(37.5,42.5)(40,40)
\psline[linewidth=0.25pt](-37.5,37.5)(-35,40)
\psline[linewidth=0.25pt](-35,35)(-32.5,37.5)
\psline[linewidth=0.25pt](-32.5,32.5)(-30,35)
\psline[linewidth=0.25pt](-30,30)(-25,35)
\psline[linewidth=0.25pt](-27.5,27.5)(-20,35)
\psline[linewidth=0.25pt](-25,25)(-17.5,32.5)
\psline[linewidth=0.25pt](-22.5,22.5)(-15,30)
\psline[linewidth=0.25pt](-20,20)(-12.5,27.5)
\psline[linewidth=0.25pt](-17.5,17.5)(-2.5,32.5)
\psline[linewidth=0.25pt](-15,15)(0,30)
\psline[linewidth=0.25pt](-12.5,12.5)(2.5,27.5)
\psline[linewidth=0.25pt](-10,10)(7.5,27.5)
\psline[linewidth=0.25pt](-7.5,7.5)(10,25)
\psline[linewidth=0.25pt](-5,5)(22.5,32.5)
\psline[linewidth=0.25pt](-2.5,2.5)(27.5,32.5)
\psline[linewidth=0.25pt](2.5,2.5)(-30,35)
\psline[linewidth=0.25pt](5,5)(-25,35)
\psline[linewidth=0.25pt](7.5,7.5)(-12.5,27.5)
\psline[linewidth=0.25pt](10,10)(-10,30)
\psline[linewidth=0.25pt](12.5,12.5)(-5,30)
\psline[linewidth=0.25pt](15,15)(2.5,27.5)
\psline[linewidth=0.25pt](17.5,17.5)(10,25)
\psline[linewidth=0.25pt](20,20)(12.5,27.5)
\psline[linewidth=0.25pt](22.5,22.5)(15,30)
\psline[linewidth=0.25pt](25,25)(17.5,32.5)
\psline[linewidth=0.25pt](27.5,27.5)(20,35)
\psline[linewidth=0.25pt](30,30)(27.5,32.5)
\psline[linewidth=0.25pt](32.5,32.5)(30,35)
\psline[linewidth=0.25pt](35,35)(32.5,37.5)
\psline[linewidth=0.25pt](37.5,37.5)(35,40)
\psline[border=1pt, bordercolor=white]{->}(0,0)(0,50)
\psline{->}(-50,0)(50,0)
\psline[linecolor=blue]{-*}(38.75,41.25)(38.75,0)
\psline[linecolor=blue]{-*}(26.25,33.75)(26.25,0)
\psline[linecolor=blue]{-*}(21.25,33.75)(21.25,0)
\psline[linecolor=blue]{-*}(6.25,28.75)(6.25,0)
\psline[linecolor=blue]{-*}(8.75,26.25)(8.75,0)
\psline[linecolor=blue]{-*}(1.25,28.75)(1.25,0)
\psline[linecolor=red]{-*}(-38.75,41.25)(-38.75,0)
\psline[linecolor=red]{-*}(-28.75,36.25)(-28.75,0)
\psline[linecolor=red]{-*}(-23.75,36.25)(-23.75,0)
\psline[linecolor=red]{-*}(-11.25,28.75)(-11.25,0)
\psline[linecolor=red]{-*}(-8.75,31.25)(-8.75,0)
\psline[linecolor=red]{-*}(-3.75,31.25)(-3.75,0)
\psline[linewidth=1pt](37.5,42.5)(40,40)
\psline[linewidth=1pt](25,35)(27.5,32.5)
\psline[linewidth=1pt](20,35)(22.5,32.5)
\psline[linewidth=1pt](5,30)(10,25)
\psline[linewidth=1pt](0,30)(2.5,27.5)
\psline[linewidth=1pt](-40,40)(-37.5,42.5)
\psline[linewidth=1pt](-30,35)(-27.5,37.5)
\psline[linewidth=1pt](-25,35)(-22.5,37.5)
\psline[linewidth=1pt](-12.5,27.5)(-7.5,32.5)
\psline[linewidth=1pt](-5,30)(-2.5,32.5)
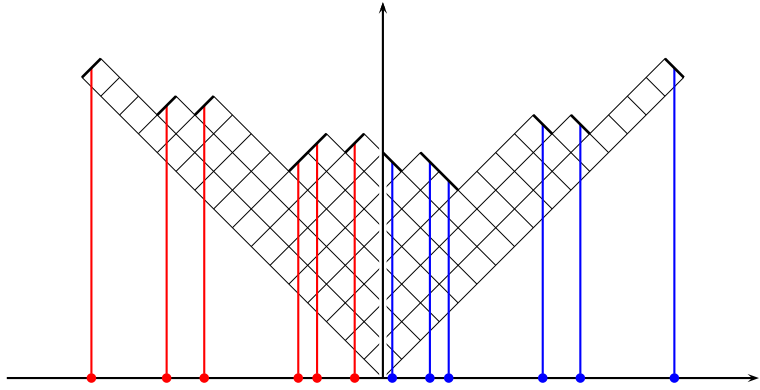
\endpspicture}{Modified Frobenius coordinates of a Young diagram drawn in Russian style. For instance, $a_1^*$ is here equal to $31/2$, and $-b_2^*=-23/2$.\label{killinginthename}}\bigskip

We denote by $\obs$ the algebra generated by the functions $(p_k)_{k \geq 1}$; it is the algebra of polynomial functions on Young diagrams introduced by S. Kerov and G. Olshanski in \cite{KO94}, but we use here the notations of \cite{IO02}. It can be shown that the $p_k$ are algebraically independent over $\C$, and as a consequence, the $p_\rho = p_{\rho_1} p_{\rho_2}\cdots p_{\rho_r}$ (where $\rho$ runs over all partitions) form a linear basis of $\obs$. The algebra of observables of diagrams has therefore a canonical grada\-tion given by:
$$\deg p_\rho=|\rho|=\rho_1+\rho_2+\cdots+\rho_r.$$
\medskip

\subsection{Normalized characters}
Another graded algebraic basis may be considered: the normalized cha\-racters $(\varSigma_k)_{k \geq 1}$. If $\rho$ and $\lambda$ are two partitions of respective sizes $k$ and $n$, we define the following function on $\young$:
$$\varSigma_\rho(\lambda)=\begin{cases}
	                  n^{\downarrow k}\,\chi^\lambda(\sigma_{\rho\,1^{n-k}})&\text{if }k\leq n\\
                          0 & \text{if } k>n,
                          \end{cases}
$$
where $n^{\downarrow k}=n(n-1)\cdots(n-k+1)$ is the $k$-th falling factorial of $n$, the notation $\chi^\lambda$ stands for the normalized character of the symmetric group, the partition $\rho'=\rho\,1^{n-k}$ is the partition $\rho$ completed with parts of size $1$, and $\sigma_{\rho'}$ is any permutation in $\sym_n$ with cycle decomposition of type $\rho'$.\bigskip

In the case where $\lambda$ has only one part (character values on cycles), by using the Frobenius formula $\varsigma^\lambda(\rho)=\scal{s_\lambda}{p_\rho}$ (\emph{cf.} \cite[\S1.7]{Mac95}, and in particular the corollary 1.7.7), it can be shown that  $\varSigma_k$ is a polynomial in the $(p_l)_{l \leq k}$. To compute these polynomials, one can use the following Proposition (see \cite[Propositions 3.3 and 3.4, 4.4 and 4.5]{IO02}):
\begin{proposition}\label{prop:sigma_pk} If $k\geq 1$, the observable $\varSigma_k$ can be expressed in terms of power sums {\it via} the following formula
\begin{equation} \label{eq:p_fct_sigma}
\varSigma_k=[t^{k+1}] \left\{-\frac{1}{k}\,\prod_{j=1}^k (1-(j-1/2)t)\,\cdot\,\exp\left(\sum_{j=1}^\infty \frac{p_j\,t^j}{j}\,(1-(1-kt)^{-j}) \right)\right\}.
\end{equation}
\noindent As a consequence, the top homogeneous component of $\varSigma_k$ is $p_k$, and $(\varSigma_k)_{k \geq 1}$ is an algebraic basis of the algebra $\obs$.
\end{proposition}
\begin{examples}
$\varSigma_{1}=p_{1}$, $\varSigma_{2}=p_{2}$,
$\varSigma_{3}=p_{3}-\frac{3}{2} p_{11}+\frac{5}{4}p_{1}$,
$\varSigma_{4}=p_{4}-4p_{21}+\frac{11}{2}p_{2}$.
\end{examples}\bigskip

\subsection{Normalized conjugacy classes in symmetric group algebras}
If $G$ is a group and $G^\vee$ a set indexing its irreducible representations,
one has a ring isomorphism, called abstract Fourier transform, between the center
of the group algebra and the algebra of functions on $G^\vee$:
\begin{align*}
    Z(\C G) &\stackrel{\sim}{\rightarrow} \Func(G^\vee,\C) \\
    x &\mapsto \big( \lambda \mapsto \chi^\lambda(x) \big)
\end{align*}
It is indeed a ring isomorphism because the element in $Z(\C G)$
acts, by Schur lemma, as an homothetic transformation on each
irreducible representation.\bigskip

In the case $G=\sym_n$, this means that a function on Young diagrams with $n$
boxes corresponds to an element of the center of the group algebra $\C\sym_n$.
Thus a function $f$ on all Young diagrams can be seen as sequences
of elements $(f_n)_{n \geq 1}$ with $f_n \in Z(\C\sym_n)$ for all $n$.
In the case of the normalized characters $\varSigma_\rho$,
this sequence has a very nice expression, see \cite{Sniady06}:
\begin{align}
\left(\varSigma_{\rho_1,\rho_2,\ldots,\rho_r}\right)_n & = \sum_{\substack{a_{(i,j)} \in \lle 1,n \rre \\ a_{i,j}\neq a_{(i,j)'}}} (a_{11},a_{12},\ldots,a_{1\rho_1})\,(a_{21},a_{22},\ldots,a_{2\rho_2})\,\cdots\,(a_{r1},a_{r2},\ldots,a_{r\rho_r}) \nonumber\\
\label{eq:sigma_alg} &= \sum_{\substack{\forall i,\,\,A_i \in \Arr(n,\,\rho_i)\\ \forall i\neq j,\,\,A_i\cap A_j=\emptyset}} C(A_1)\,C(A_2)\,\cdots \,C(A_r),
\end{align}
where $\Arr(n,k)$ is the set of $k$-arrangements of $\lle 1,n \rre$
and $C(A)$ is the cycle corresponding to a sequence $A$ of distinct numbers.
It is easy to check that these elements verify
$$\forall \lambda \in \young, \ \ \varSigma_\rho(\lambda)
=\chi^\lambda\big((\varSigma_\rho)_{|\lambda|}\big).$$
Hence, the element in equation (\ref{eq:sigma_alg}) corresponds to the
function $\varSigma_\rho$ on Young diagrams {\it via} the abstract Fourier transform,
so the notation is coherent. Moreover, formula \eqref{eq:sigma_alg}
allows us to do computations for all $n$ simultaneously and in this case,
we will forget the index $n$ (for example in next paragraph).\bigskip

\subsection{Products of normalized characters}
 As the abstract Fourier transform is a ring homomorphism, 
the point-wise multiplication of functions over Young diagrams of size $n$
 corresponds to the product in the center of the symmetric group algebra.
 But, in $Z(\C\sym_n)$, one has
$$\varSigma_\mu \cdot \varSigma_\nu = \sum_{\substack{\forall i,\,\,A_i \in \Arr(n,\,\mu_i)\\ \forall i\neq j,\,\,A_i\cap A_j=\emptyset}} \ \sum_{ \substack{\forall i,\,\,B_i \in \Arr(n,\,\nu_i)\\ \forall i\neq j,\,\,B_i\cap B_j=\emptyset}} C(A_1)\,C(A_2)\,\cdots \,C(A_r)\, C(B_1)\,C(B_2)\,\cdots \,C(B_r).$$
As the $A_i$'s and the $B_j$'s may have elements in common, this is not equal to $\varSigma_{\mu \cup \nu}$. Let us split this double summation over $a_{i,j}$ and $b_{k,l}$ in the following way. To each indexing sequence $a_{i,j}, b_{k,l}$, we associate the partial matching between the sets $I_A=(i,j)_{1 \leq j \leq \mu_i}$ (indices of the $a$'s) and $I_B=(k,l)_{1 \leq k \leq \nu_l}$, which links $(i,j)$ and $(k,l)$ if and only if $a_{i,j}=b_{k,l}$. When we sum the quantity 
\begin{equation}\label{eq:prod_cycles}
 C(A_1)\,C(A_2)\,\cdots \,C(A_r)\, C(B_1)\,C(B_2)\,\cdots \,C(B_r)
\end{equation}
over indexing sequences with a given partial matching, we obtain $\Sigma_\rho$ for some $\rho$ depending on the matching. Indeed, if we know which of the $a_{i,j}$'s are equal to which of the $b_{k,l}$'s, the product \eqref{eq:prod_cycles} can be rewritten as a product of cycles with disjoint supports in a way which does not depend of the values of the $a$'s and of the $b$'s. Finally,
$$\varSigma_\mu \cdot \varSigma_\nu = \sum \varSigma_{\rho(M)},$$
where the sum runs over partial matchings  $M$ of the set of indices of the $a$'s with the set of indices of the $b$'s.\bigskip

\begin{example}
Let us compute $\varSigma_{2} \cdot \varSigma_{3}$. The set of indices of the $a$'s (resp. the $b$'s) is $\{1,2\}$ (resp. $\{1,2,3\}$).
\begin{itemize}
 \item If $A$ and $B$ have no common element, the product $C(A) \, C(B)$ is a 3-cycle multiplied by a 2-cycle with disjoint supports and the corresponding sum is $\varSigma_{3,2}$.
 \item If $a_1=b_1$, but all other elements are different, then $(a_1\ a_2\ a_3)\cdot (a_1\ b_2) = (a_1\ b_2\ a_2\ a_3)$. When we sum on all possible values of $a_1$, $b_2$, $a_2$, $a_3$ (all different), one obtains $\varSigma_4$.

The result is the same for any of the $6$ matchings of size $1$.
 \item If $a_1=b_1$ and $a_2=b_2$, then $(a_1,a_2,a_3)\cdot (a_1,a_2) = (a_1,a_3)(a_2)$. When we sum on all possible values of $a_1$, $a_2$, $a_3$, we obtain $\varSigma_{2,1}$ (note that here it is important to remember that we have to sum over indices $a_2$ even if its value does not change the permutation).

Once again, the result is the same for any of the $6$ matchings of size $2$.
\end{itemize}

Thus, $\varSigma_{2} \cdot \varSigma_{3} = \varSigma_{2,3} + 6 \varSigma_4 + 6 \varSigma_{2,1}$.
\end{example}\medskip

\noindent In general, one has the following properties:
\begin{enumerate}
 \item The size of $\rho(M)$ is simply $|\mu|+|\nu|-|M|$ where $|M|$ is the size of the matching.
 \item The term corresponding to the empty matching is $\varSigma_{\nu \cup \mu}$. Thus, $$\varSigma_\nu \cdot \varSigma_\mu = \varSigma_{\nu \cup \mu} + \sum_{|\rho| < |\mu|+ |\nu|} c^\rho_{\nu,\mu} \,\varSigma_\rho.$$ By induction, one obtains that, if $k_1 \geq k_2 \geq \ldots \geq k_r$:
$$\varSigma_{k_1} \cdot \varSigma_{k_2}\cdots \varSigma_{k_r} = \varSigma_{k_1,\ldots,k_r} + \sum_{|\rho| < k_1 + \cdots + k_r} c^\rho_{k_1,\ldots,k_r}\, \varSigma_\rho.$$
This triangular relation can be inverted, so the $\varSigma_\rho$ are in the algebra $\obs$ and span it linearly. Moreover, $\varSigma_\rho$ is a non-homogeneous element of degree $|\rho|$ and its homogeneous component of maximal degree is $p_\rho$.
 \item If $\mu=(l)$ and $\nu=(m)$, all matchings of size $1$ give a term $\varSigma_{l+m-1}$. As there are $ml$ such matchings, one has $$\varSigma_l \,\varSigma_m = \varSigma_{l,m} +ml\, \varSigma_{l+m-1}+(\text{terms of degree at most }l+m-2).$$
\end{enumerate}\bigskip

\subsection{Another product}\label{subsect:bullet_product}

In this paragraph, we recall the definition of the disjoint product $\bullet$ on $\obs$ (see \cite{Sniady06}), which will be useful in section \ref{sect:second_order}. It is enough to give the value of the product of two elements of the basis $\varSigma_\rho$:
$$\varSigma_\mu \bullet \varSigma_\nu = \varSigma_{\mu \cup \nu}.$$
If we look at the description of the classical product on the $\varSigma_\mu$ given at the previous paragraph, taking the product $\bullet$ consists in keeping in the result the terms indexed by families of disjoint indices. This remark can be generalized to expression of the kind:
$$X_{i_1,\dots,i_r}^{\mathcal{G}} = \sum_{ \substack{\forall j,\,\,A_j \in \Arr(n,i_j) \\  G(a_{i,j}) \in \mathcal{G}}} \!\!\!\!\!\!\!\! C(A_1) \cdots C(A_r),$$
where $G(a_{i,j})$ is, as in the previous paragraph, the graph describing which of the $a_{i,j}$'s are equal and $\mathcal{G}$ a collection of \emph{admissible} graphs. The following lemma will be useful in subsection \ref{subsect:bound_cumulants}:
\begin{lemma}\label{lem:bullet}
The disjoint product of two such expressions is given by:
$$X_{i_1,\dots,i_r}^{\mathcal{G}} \bullet X_{i'_1,\dots,i'_r}^{\mathcal{G}'} = X_{i_1,\dots,i_r,i'_1,\dots,i'_r}^{\mathcal{G} \sqcup \mathcal{G}'},$$
where $\mathcal{G} \sqcup \mathcal{G}'$ is the collection of all graphs obtained by the disjoint union of a graph of $\mathcal{G}$ and a graph of $\mathcal{G}'$.
In other words, the disjoint product of two expressions as above is obtained by keeping in the usual product only terms such that the two families of indices are disjoint.
\end{lemma}
\begin{proof}
As in the previous paragraph, the summation over sequences of indices corresponding to a given graph gives a term $\varSigma_\rho$. Therefore:
$$X_{i_1,\dots,i_r}^{\mathcal{G}} = \sum_{G \in \mathcal{G}} \varSigma_{\rho(G)}$$
and we have
\begin{align*}
X_{i_1,\dots,i_r}^{\mathcal{G}} \bullet X_{i'_1,\dots,i'_r}^{\mathcal{G}'} & = \sum_{G \in \mathcal{G}} \sum_{G' \in \mathcal{G}'} \varSigma_{\rho(G),\rho(G')};\\
&=X_{i_1,\dots,i_r,i'_1,\dots,i'_r}^{\mathcal{G} \sqcup \mathcal{G}'}.\qedhere
\end{align*}
\end{proof}\bigskip

\subsection{$q$-characters}

In this paragraph, we present a quantization of the $\varSigma$-basis of $\obs$
corresponding to renormalized $q$-characters of the Iwahori-Hecke algebras.
The theory of character of Iwahori-Hecke algebra is a little more involved
than the one of symmetric groups.
Indeed, given two permutations $\sigma$ and $\tau$ in the same conjugacy class of $\sym_n$,
the $T$-basis elements $T_\sigma$ and $T_\tau$ can have different traces in a $\IH_{n,q}$-module:
for instance, $(1,2)$ and $(1,3)$ are conjugate in $\sym_3$, but $\chi(T_{(1,2)})=q$
and $\chi(T_{(1,3)})=q^3$ if we consider the one-dimensional index representation.
However, if $\sigma$ and $\tau$ are of minimal Coxeter length in their conjugacy classes,
then one has $\chi(T_{\sigma})=\chi(T_{\tau})$ for any character $\chi$.
Moreover, the values of a character $\chi$ on these particular elements determine $\chi$,
see  \cite[Propositions 5.1 and 5.2]{Ram91} or \cite[\S8.2]{GP00} for the general case
of an Iwahori-Hecke algebra of a Coxeter group. Consequently, if
$$\sigma_\rho=(1,2,\ldots,\rho_1)\,(\rho_1+1,\rho_1+2,\ldots,\rho_1+\rho_2)\,
\cdots\,(\rho_1+\cdots+\rho_{r-1}+1,\ldots,\rho_1+\cdots+\rho_r),$$
then $\varsigma^{\lambda}(q)$ is entirely determined by the values
$\varsigma^{\lambda}(\rho,q)=\varsigma^{\lambda}(T_{\sigma_{\rho}},q)$, 
where $\rho$ runs over $\young_{n}$.
As in the case $q=1$, we denote by $\chi^\lambda(\rho,q)$ the rescaled
character value (divided by the dimension of the module).\bigskip

In \cite{Ram91}, A. Ram uses the Schur-Weyl decomposition of $(\C^m)^{\otimes n}$
into irreducible $(U_q(\mathfrak{sl}(m,\C)),\IH_{n,q})$-bimodules
to compute explicitly the characters of the Iwahori-Hecke algebra $\IH_{n,q}$ 
--- the reader can also consult \cite{RRW96} and \cite{RR97}.
These computations lead to a generalization of the usual Frobenius-Schur formula. \bigskip

To state it, we need to introduce a few notations about symmetric functions (we try
to be consistent with the one in \cite{Mac95}): the notations $h_\lambda$, $m_\lambda$,
$s_\lambda$ and $p_\lambda$ stand respectively for the complete, monomial, Schur and
power-sum basis of the symmetric function ring. We will also use the Hall scalar product,
for which complete function basis is dual to the monomial one, while Schur functions are
orthonormal and power sums orthogonal.\bigskip

Furthermore, let us define a $q$-analog of the power sum $p_k(X)$.
Let $\widetilde{q}_k(X,q)$ be the symmetric function defined by the formulas
\begin{align*}
\widetilde{q_k}(X,q)&=\frac{q^k}{q-1}\, h_k(X(1-q^{-1})) \quad (\lambda\text{-ring notations, }X\text{ is an alphabet and }q^{-1}\text{ is a variable}),\\
&=\frac{q^k}{q-1}\,\sum_{\mu \in \young_k} (1-q^{-1})^{\ell(\mu)} \,m_\mu(X)=\frac{1}{q-1} \,\sum_{\mu \in \young_k} (q^\mu-1)(z_\mu)^{-1}\,p_\mu(X),
\end{align*}
where $q^\mu-1$ means $\prod_{i=1}^r \,(q^{\mu_i}-1)$, and where $z_\mu$ equals $\prod_{i \geq 1} i^{m_i}\,m_i!\,$ if $\mu$ writes multiplicatively as $1^{m_1}2^{m_2}\ldots s^{m_s}$. This is a slightly modified version of the Hall-Littlewood polynomial $P_k(X,q^{-1})$ described in \cite[Chapter 3]{Mac95}. Let us first recall Frobenius formula for usual character values of symmetric group:
$$\forall \rho \in \young_n,\,\,\,p_\rho(X)=\sum_{\lambda \in \young_n} \varsigma^{\lambda}(\rho)\,s_\lambda(X).$$
Here is the statement of Ram's generalization of Frobenius-Schur formula:
\begin{proposition}[$q$-Frobenius-Schur formula]\label{quantumcross}
If $\widetilde{q}_\rho=\widetilde{q}_{\rho_1}\cdots \widetilde{q}_{\rho_r}$, then the following generalized Frobenius-Schur formula holds:
$$\forall \rho \in \young_n,\,\,\,\widetilde{q}_\rho(X,q)=\sum_{\lambda \in \young_n} \varsigma^{\lambda}(\rho,q)\,s_\lambda(X).$$
\end{proposition}\bigskip

Many combinatoric rules follow from Proposition \ref{quantumcross}, \emph{e.g.}, a quantization of the classical Murnaghan-Nakayama rule, see \cite[\S2]{RR97}. For the asymptotic analysis of the $q$-Plancherel measures, we shall use a triangular relation between $q$-characters and standard characters:
\begin{align*}
\varsigma^{\lambda}(\rho,q)&=\scal{\widetilde{q}_\rho(X,q)}{s_\lambda(X)} =\frac{q^n}{(q-1)^{\ell(\rho)}}\,\scal{h_\rho(X-q^{-1}X)}{s_\lambda(X)}\\
&=\frac{q^n}{(q-1)^{\ell(\rho)}}\,\sum_{\nu\in \young_n} \frac{\scal{h_\rho(X-q^{-1}X)}{p_\nu(X-q^{-1}X)}\,\scal{p_\nu(X-q^{-1}X)}{s_\lambda(X)}}{\scal{p_\nu}{p_\nu}}\\
&=\frac{1}{(q-1)^{\ell(\rho)}}\,\sum_{\nu \in \young_n}(q^\nu-1)\,\frac{\scal{h_\rho}{p_\nu}}{z_\nu}\,\scal{p_\nu}{s_\lambda}=\frac{1}{(q-1)^{\ell(\rho)}}\,\sum_{\nu \in \young_n}\frac{q^\nu-1}{z_\nu}\scal{h_\rho}{p_\nu}\,\varsigma^\lambda(\nu)
\end{align*}
As for $q=1$, we define the following function on the set $\young$ of all Young diagrams:
$$\varSigma_{\rho,q}(\lambda)= \left\{\begin{array}{ll} |\lambda|^{\downarrow |\rho|}\,\chi^\lambda(\rho\,1^{|\lambda|-|\rho|},q) & \text{ if }|\rho| \leq |\lambda| ; \\ 0 & \text{ else.} \end{array}\right.$$
If $|\rho|=|\lambda|=n$, the previous formula can be rewritten:
\begin{equation}\label{eq:sigmaq_fct_sigma}
(q-1)^{\ell(\rho)}\,\varSigma_{\rho,q}(\lambda)=\sum_{\nu \in \young_n} \frac{q^\nu-1}{z_\nu}\scal{p_\nu}{h_\rho}\,\varSigma_\nu(\lambda)
\end{equation}
This is of course also true for $|\rho| > |\lambda|$. We will show that it still holds for $k=|\rho| < |\lambda|=n$. In this case, the scalar product $\langle h_{\rho\,1^{n-k}}| p_\nu \rangle$ in the sum above is given by:
\begin{multline*}
\langle h_{\rho\,1^{n-k}}| p_\nu \rangle=\scal{h_\rho(p_1)^{n-k}}{p_\nu}=\sum_{\pi_1\in \young_{\rho_1},\ldots,\pi_m\in \young_{\rho_m}} (z_{\pi_1}\cdots z_{\pi_m})^{-1} \,\scal{p_\pi(p_1)^{n-k}}{p_\nu} \\
=\sum_{\pi_1,\ldots,\pi_m} z_{\pi1^{n-k}} \,(z_{\pi_1}\cdots z_{\pi_m})^{-1}\,\mathbb{1}_{(\nu=\pi1^{n-k})}
=(z_\nu\,\mathbb{1}_{(\nu=\nu_k1^{n-k})})\,\sum_{\pi_1,\ldots,\pi_m} (z_{\pi_1}\cdots z_{\pi_m})^{-1}\,\mathbb{1}_{(\nu_k=\pi)}\\
=(z_\nu /z_{\nu_k})\,\mathbb{1}_{(\nu=\nu_k1^{n-k})}\,\sum_{\pi_1,\ldots,\pi_m} (z_{\pi_1}\cdots z_{\pi_m})^{-1}\,\scal{p_\pi}{p_{\nu_k}}=(z_\nu /z_{\nu_k})\,\mathbb{1}_{(\nu=\nu_k1^{n-k})}\,\scal{h_\rho}{p_{\nu_k}}
\end{multline*}
As a consequence, formula (\ref{eq:sigmaq_fct_sigma}) can be seen as an equality of functions on $\young$. Note that this relation between the symbols $\varSigma_{\rho,q}$ and $\varSigma_\nu$ is triangular, because $\scal{p_\nu}{h_\rho}=0$ unless $\nu$ is finer than  $\rho$. Moreover, as $(m_\rho)_{\rho \in \young}$ is the dual basis of $(h_\rho)_{\rho \in \young}$, the inverse relation is easily obtained. We conclude that:
\begin{proposition}[Quantization of the algebra of observables]
 The symbols $(\varSigma_{\rho,q})_{\rho \in \young}$ form a $\C(q)$-basis of $\obs(q)=\obs \otimes_{\C}\C(q)$ and yield a quantization of the basis $(\varSigma_\rho)_{\rho \in \young}$ of $\obs$. Moreover, the transition matrices between the $\varSigma_{\rho,q}$ and the $\varSigma_\rho$ are triangular, given by:
\begin{align}
(q-1)^{\ell(\rho)}\,\varSigma_{\rho,q}(\lambda)&=\sum_{\nu \in \young_k} \frac{q^\nu-1}{z_\nu}\scal{p_\nu}{h_\rho}\,\varSigma_\nu(\lambda) ;\\
(q^\rho-1)\,\varSigma_\rho(\lambda)&=\sum_{\nu \in \young_k} (q-1)^{\ell(\nu)}\scal{m_\nu}{p_\rho}\,\varSigma_{\nu,q}(\lambda).
\end{align}

\end{proposition}\bigskip

\section{First order asymptotics}\label{sect:first_order}
This section is devoted to the proof of Theorem \ref{th:first_order_as_rows}. Let us fix a real $q \in \,]0,1[$. The first step of the proof is the convergence of power sums of Frobenius coordinates.

\subsection{Convergence of power sums}
The $q$-Plancherel measure is a measure on Young diagrams, so the functions on Young diagrams and in particular the elements of $\obs$ are in this context random variables. In particular we will denote by $\esper$ their expectation. Thanks to the characterization of the $q$-Plancherel measure with the normalized trace (paragraph \ref{subsect:qPlancherel_trace}), one has immediately:
$$\esper[\varSigma_{\rho,q}] = \begin{cases}
                                n^{\downarrow k} & \text{if }\rho = 1^k ;\\
				0 & \text{else.}
                               \end{cases}$$
By using the triangular relation between the $\varSigma_{\rho}$ and the $\varSigma_{\nu,q}$, we obtain:
\begin{align}\esper[\varSigma_{\rho}] &= \frac{1}{q^\rho-1}\,\sum_{\nu \in \young_k} (q-1)^{\ell(\nu)}\,\scal{m_\nu}{p_\rho}\,\esper[\varSigma_{\nu,q}]=\frac{(q-1)^{|\rho|}}{q^\rho-1}\,n^{\downarrow |\rho|}\,\scal{m_{1^{|\nu|}}}{p_\rho} \nonumber\\
\label{eq:esp_sigma} &=\frac{(1-q)^{|\rho|}}{1-q^\rho}\,n^{\downarrow|\rho|}
\end{align}\medskip

\noindent This gives in particular the asymptotic order of magnitude:
$$\esper[\varSigma_{\rho}] = O\big(n^{|\rho|}\big) = O\big(n^{\deg(\varSigma_{\rho})}\big)$$
In fact, this is true for any measure on Young diagrams, but is not always precise enough (for example, to study the usual Plancherel measure, one has to use the fact that it is $O(n^{(|\rho|+\ell(\rho))/2})$). In this context, the exact computation above shows that we could not find a better bound.
As $(\varSigma_{\rho})_{\rho \in \young}$ is a linear basis of $\obs$, the previous estimate remains true for any element of $\obs$.
\begin{lemma}\label{lem:major_Ex_deg}
 For any $x \in \obs$, one has
$\esper[x] = O(n^{\deg(x)})$.
\end{lemma}
\noindent This lemma allows us to forget terms of lower degree when we want to look at the asymptotic behavior of the expectation of some observable. In particular, one can obtain the following result for the power sums of the Frobenius coordinates:

\begin{lemma}[Convergence of power sums]\label{lem:cv_powersums}
Under the $q$-Plancherel measures, we have convergence in probability of the renormalized power sums $p_k(\lambda)$:
$$\forall k \geq 1,\,\,\frac{p_k(\lambda)}{|\lambda|^k} \longrightarrow_{M_{n,q}} \frac{(1-q)^k}{1-q^k}.$$
\end{lemma}
\begin{proof}
An immediate consequence of equation (\ref{eq:esp_sigma}) is that $\esper[\varSigma_\rho]/n^{|\rho|} \to(1-q)^{|\rho|}/(1-q^\rho)$ when $n$ goes to infinity. But $\varSigma_\rho-p_\rho$ is an observable of degree less than $|\rho|-1$, so the same holds for $p_\rho$:
\begin{equation}\label{EqEquivPRho}
\lim_{n \to \infty} \frac{\esper[p_\rho]}{n^{|\rho|}} = \frac{(1-q)^{|\rho|}}{1-q^\rho}
\end{equation}
Fix $\eps>0$ and $k \geq 1$. The Bienaym\'e-Chebyshev inequality ensures that
\begin{align*}M_{n,q}\left[\left|\frac{p_k}{n^k} - \frac{(1-q)^k}{1-q^k}\right|\geq \eps\right] &\leq\frac{1}{\eps^2}\,\esper\left[\left(\frac{p_k}{n^k} - \frac{(1-q)^k}{1-q^k}\right)^2\right]\\
&\leq\frac{1}{\eps^2}\,\left(\frac{\esper[p_{k,k}]}{n^{2k}} - 2 \,\frac{\esper[p_k]}{n^k}\,\frac{(1-q)^k}{1-q^k} +\frac{(1-q)^{2k}}{1-q^{k,k}}\right).
\end{align*}
But the right-hand side tends to $0$ (which is easily seen using equation \eqref{EqEquivPRho} for $\rho=(k)$ and $\rho=(k,k)$) and the proof is over.
\end{proof} \bigskip

\begin{corollary}[Convergence of the characters and $q$-characters]\label{corol:concentration_q-characters}
Under the $q$-Plancherel measures, in probability,o ne has:
\begin{align*}
 &\forall \rho \in \young, \,\, \frac{\varSigma_\rho}{n^{|\rho|}}  \longrightarrow_{M_{n,q}} \frac{(1-q)^{|\rho|}}{1-q^\rho}\quad   \text{  i.e. } \chi^\lambda(\rho)  \longrightarrow_{M_{n,q}} \frac{(1-q)^{|\rho|}}{1-q^\rho};\\ 
& \forall \rho \in \young,  \,\, \frac{\varSigma_{\rho,q}}{n^{|\rho|}}  \longrightarrow_{M_{n,q}} \delta_{\rho=1^{|\rho|}}\quad \quad  \text{ i.e. } \chi^\lambda(\rho,q)  \longrightarrow_{M_{n,q}}\delta_{\rho=1^{|\rho|}}.
\end{align*}
\end{corollary}

\begin{proof}
It is a direct consequence of Proposition \ref{prop:sigma_pk} and Lemma \ref{lem:cv_powersums}.
\end{proof}\bigskip

\subsection{Convergence of row lengths} \label{subsect:convergence_rows}
 Note that the limit values of the rescaled $p_k$'s are the power sums of the sequence $((1-q),(1-q)q,(1-q)q^2,\ldots)$. Indeed,
$$\sum_{i=1}^\infty (1-q)^k\, q^{(i-1)k}=\frac{(1-q)^k}{1-q^k}.$$
In this paragraph, we show that the convergence of power sums implies the convergence of the sequence. This corresponds to the end of the proof of Theorem \ref{th:first_order_as_rows}. \bigskip

The main idea is to consider, given a Young diagram $\lambda$ of size $n$, the probability measure $X_\lambda$ defined by
$$X_\lambda=\sum_{i=1}^d (a_i^*(\lambda)/n)\,\delta_{(a_i^*(\lambda)/n)}+(b_i^*(\lambda)/n)\,\delta_{(-b_i^*(\lambda)/n)}.$$
The $k$-th moment of $X_\lambda$ is exactly $p_{k+1}(\lambda)/n^{k+1}$. So, thanks Lemma \ref{lem:cv_powersums}, under $q$-Plancherel measure, the moments of $X_\lambda$, which can be seen as random variables, converge in probability towards the moments of the measure
$$ X_{\infty;q}=\sum_{i=1}^\infty (1-q)q^{i-1}\,\delta_{(1-q)q^{i-1}}. $$
\begin{lemma}\label{lem:cv_fct_repartition}
As usual, we consider a sequence of random diagrams $\lambda_n \in \young_n$ taken with $q$-Plancherel measure. At any point $x \neq (1-q)\,q^{i-1}$, one has convergence in probability of the repartition functions, that is to say:
$$\sum_{\substack{y\leq x\\ ny \in A_\lambda+\overline{B}_\lambda}} \!\!\!\!|y|\quad \longrightarrow_{M_{n,q}} \sum_{(1-q)q^{i-1} \leq x} \!\!\!\!(1-q)q^{i-1}.$$
\end{lemma}
 
\begin{proof}
In the following, we consider $\lambda \mapsto X_\lambda$ as a \emph{measure-valued random variable} $X_{n,q}$; its law  is by definition the image of $M_{n,q}$ by $\lambda \mapsto X_\lambda$. The set $\mathscr{M}([-1,1])$ of probability measures on $[-1,1]$, where $X_{n,q}$ takes its values, is included in the dual of $\mathscr{C}([-1,1])$ and is compact and metrizable with respect to the $*$-weak topology\footnote{It is a particular case of the Prohorov theorem: if $E$ is a separable and complete metric space (in other words, a polish space), then the same holds for the set of probability measures $\mathscr{M}(E)$, see \cite[Chapter 1]{Bil69}. As a consequence, convergence in probability makes sense for random probability measures on a closed subset of the real line.}. Indeed, if $(f_n)_{n \in \N}$ is a sequence of continuous functions dense in $\mathscr{C}([-1,1])$, then
$$d(m_1,m_2)=\sum_{n \in \N} \frac{1}{2^n} \,\max(1,|m_1(f_n)-m_2(f_n)|)$$
is a distance compatible with the topology of $\mathscr{M}([-1,1])$. We choose $\{f_n\}_{n \in \N}=\mathbb{Q}[X]$ (dense by Stone-Weierstrass theorem), thus the convergence of moments implies that $X_{n,q}$ converges in probability in $(\mathscr{M}([-1,1]),d)$ towards $X_{\infty;q}$.\bigskip

\noindent But, for Radon measures on the real line, the weak convergence $X_{n,q} \rightharpoonup X_{\infty;q}$ is equivalent to the convergence of repartition functions $F_{X_{n,q}}(x) \to F_{X_{\infty;q}}(x)$ at points $x \in [-1,1]$ where $F_{X_{\infty;q}}$ is continuous --- it is a part of the Portmanteau theorem (\cite[Chapter1]{Bil69}).
\end{proof}\bigskip

\begin{proof}[End of proof of Theorem \ref{th:first_order_as_rows}]
Now, let us take $x=1-q+\eta$ with $\eta >0$ sufficiently small. We also fix $1-q>\eps>0$. Note that $F_{X_{\infty;q}}=1$ on an open set containing $x$. Consequently, for $n$ big enough and outside a set of arbitrarily small probability, all rows of $\lambda$ are smaller than $nx$, because
$$F_{X_\lambda}(x)=\sum_{\substack{y \leq x \\ ny \in A+\overline{B}}} \!\!\!\!|y|\,\,\, \geq F_{X_{\infty;q}}(x)-\eps = 1-\eps > q,$$
and this could not be if a row were bigger than $n(1-q+\eta)$. Indeed, we would then have
$$F_{X_\lambda}(1) \geq F_{X_{n,q}}(x)+1-q+\eta > 1+\eta>1.$$
Hence, for $n$ big enough, $a_1^*(\lambda)/n$ is smaller than $x$ outside a set of arbitrary small probability. Now, if $x'=1-q-\eta$, then $F_{X_{\infty;q}}=q$ on an open set containing $x'$, and therefore, for $n$ big enough and outside a set of arbitrarily small probability, some rows of $\lambda$ are bigger than $nx'$; otherwise, we would have $$F_{X_\lambda}(1)=F_{X_\lambda}(x')\leq F_{X_{\infty;q}}(x')+\eps= q+\eps <1.$$ We conclude that for $n$ big enough and outside a set of arbitrarily small probability, the rescaled first column $a_1^{*}(\lambda)/n$ is between $x'$ and $x$. Consequently, $a_1^*(\lambda)/n$ converges in probability towards $1-q$, and the same argument holds for the following rows (with an induction on the index $i$ of the row). Of course, $a_i^*(\lambda)/n-\lambda_i/n=O(1/n)$, so we have proved Theorem \ref{th:first_order_as_rows}. \end{proof}\bigskip

\begin{remark}\label{RemColumns}
Unfortunately, our estimates are not strong enough to prove almost surely convergence.

Let us look at the first column of a random diagram under the $q$-Plancherel measure.
 If we take $x=0$ in Lemma \ref{lem:cv_fct_repartition}, then we see that $\sum_{i=1}^d b_i^*(\lambda)/n$ goes to $0$ in probability: so, the column lengths are $o(n)$. Nevertheless, we can not determine their precise order of magnitude.\bigskip
\end{remark}

\begin{remark}
 The same proof works for any measure for which we can compute expectation of character values as soon as (the following condition implies the convergence of variance of power sums towards $0$ in Lemma \ref{lem:cv_powersums}):
$$\esper[\Sigma_{k,k}] = \esper[\Sigma_k]^2 + o(n^{2k}).$$
This is similar to the result obtained by Biane for balanced diagrams in \cite{Bia01}.
\end{remark}
\bigskip

\subsection{Link with the representation theory of the infinite symmetric group}

This result on the first order asymptotics could alternatively have been obtained
by using general results coming from the representation theory of the infinite
symmetric group $\sym_{\infty}$, see in particular \cite{KV81,Ker03,KOV04}.
However, it seems impossible to understand the fluctuations with this approach.

Recall that the Thoma simplex 
$$\Omega=\left\{\big((\alpha_i)_{i \geq 1},(\beta_i)_{i \geq 1}\big) \,\,\big|\,\, \alpha_1\geq \alpha_2\geq \cdots \geq 0,\,\, \beta_1\geq \beta_2\geq \cdots \geq 0,\,\, \gamma=1-\sum_{i=1}^\infty (\alpha_i+\beta_i) \geq 0\right\}$$
parametrizes the normalized irreducible characters of $\sym_{\infty}$ (\cite{Tho64}).
Besides, to each point $\omega \in \Omega$, one can associate a system $(M^{\omega}_{n})_{n\in \N}$
of probability measures on the sets $\young_{n}$ of Young diagrams:
$$M^{\omega}_{n}(\lambda)=(\dim \lambda)\,s_{\lambda}(\omega),$$
 where $s_{\lambda}(\omega)$ is the generalized Schur function defined by the following specialization of the algebra of symmetric functions:
$$p_{1}(\omega)=1\ ;\qquad p_{k\geq 2}(\omega)=p_{k}(\alpha-(-\beta))=\sum_{i\geq 1}(\alpha_{i})^{k}+(-1)^{k-1}\sum_{i \geq 1}(\beta_{i})^{k}.$$
This system is called coherent because it fulfills the relation:
$$\forall n,\,\,\,\forall \lambda \in \young_{n},\,\,\,M^{\omega}_{n}(\lambda)=\sum_{\lambda \nearrow \Lambda} \frac{\dim \lambda}{\dim \Lambda}\,M^{\omega}_{n+1}(\Lambda).$$
%Then, the restriction to $\sym_{n}$ of the irreducible character $\chi^{\omega}$ of $\sym_{\infty}$ is $\sum_{\lambda \vdash n}M^{\omega}_{n}(\lambda)\,\chi^{\lambda}$, so:
%$$\chi^{\omega}(\mu)=\sum_{n}\varsigma^{\lambda}(\mu)\,s_{\lambda}(\omega)=p_{\mu}(\omega)=\prod_{k \geq 2}\left(\sum_{i \geq 1}(\alpha_{i})^{k} +(-1)^{k-1}\sum_{i \geq 1}(\beta_{i})^{k}\right)^{m_{k}(\mu)}.$$
Moreover, one can describe all coherent systems $(M^{m}_{n})_{n\in \N}$ of probability measures on the sets $\young_{n}$: given a probability measure $m \in \mathscr{M}(\Omega)$, let us define:
$$M^{m}_{n}(\lambda)=\int_{\Omega}M^{\omega}_{n}(\lambda)\,dm(\omega)=\int_{\Omega}(\dim\lambda)\,s_{\lambda}(\omega)\,dm(\omega).$$
This system is coherent and the correspondence is one-to-one.
In other terms $\Omega$ is the Martin boundary of the Young graph (\cite{KOO97}).
The set of probability measures $\mathscr{M}(\Omega)$ is also in correspondence
with normalized (but non necessarily irreducible) characters of $\sym_{\infty}$.

A deep result due to Kerov and Vershik ensures that the Thoma simplex is also in some sense the ``geometric boundary'' of the Young graph. More precisely, let us restate here Theorem 9.7.3 of \cite{KOV04} (we refer also to \cite{KV81}). To each Young diagram $\lambda \in \young_{n}$, we associate a point $\omega_{n}(\lambda)$ of the Thoma simplex by scaling its Frobenius coordinates:
$$\omega_{n}(\lambda)=\left(\frac{a_{1}^{*}(\lambda)}{n},\ldots,\frac{a_{d}^{*}(\lambda)}{n},0,\ldots ;\frac{b_{1}^{*}(\lambda)}{n},\ldots,\frac{b_{d}^{*}(\lambda)}{n},0,\ldots\right).$$
We endow $\Omega$ with the compact topology of simple convergence of coordinates.
\begin{theorem}
Let $m$ be a probability measure on $\Omega$; we denote by $(M^{m}_{n})_{n \in \N}$ the associated coherent system of measures. The push-forward $(\omega_{n})_{\star}M^{m}_{n}$ converges in law towards $m$ in $\mathscr{M}(\Omega)$.
\end{theorem}\bigskip

That said, the $q$-Plancherel measures $(M_{n,q})_{n \in \N}$ form the coherent system associated to the Dirac mass at the point 
$$\omega_{q}=\big((1-q),(1-q)q,(1-q)q^{2},\ldots;0,0,\ldots\big).$$
Indeed, for any partition $\lambda$, the generic degree $D_{\lambda}(q)$ can be written as $D_{\lambda}(q)=\{n!\}_{q}\,s_{\lambda}(\omega_{q}) $, see \cite[Chapter 1,\S3, examples 1-5]{Mac95}. So, the scaled Frobenius coordinates under $q$-Plancherel measures should converge towards the vector $(1-q,(1-q)q,(1-q)q^{2},\ldots;0,0,\dots)$, which is exactly our Theorem \ref{th:first_order_as_rows} and Remark \ref{RemColumns}.
\bigskip

\section{Second order asymptotics}\label{sect:second_order}

In this section, we investigate the second order asymptotics of the rows of diagrams under the $q$-Plancherel measure and we prove Theorem \ref{th:second_order_as_rows}. As in the previous section, the first part of the proof consists in studying asymptotics of $p_k$. In particular, we shall prove:

\begin{proposition}[Fluctuations of $p_k$]\label{prop:deviation_powersums}
Under the $q$-Plancherel measure, we have the following convergence result. If we denote the renormalized centered power sum by
$$W_{n,q,l}(\lambda)=\sqrt{n}\left(\frac{p_l(\lambda)-\esper[p_l]}{n^l}\right),$$
then, for any $k \geq 1$, the finite random vector $(W_{n,q,1},W_{n,q,2},\ldots,W_{n,q,k})$ converges in law towards a gaussian vector of covariance matrix:
$$\mathrm{cov}(W_{q,l},W_{q,m})=lm\,(1-q)^{l+m}\,\left(\frac{1}{1-q^{l+m-1,1}}-\frac{1}{1-q^{l,m}}\right).$$
\end{proposition}\bigskip 

\noindent Before we start, recall that the algebra of observables $\obs$ admits two different products: the \emph{usual product} of (commuting) observables $\cdot$ and the \emph{disjoint product} $\bullet$, \emph{cf.} section \ref{sect:observables}. We shall denote by $\obs^\bullet$ the algebra with the latter product, and by $\esper^\mathrm{id}$ the identification $\obs \to \obs^\bullet$. The proof of Proposition \ref{prop:deviation_powersums} relies on the interaction of the two products and on the expectation of character values $\esper[\varSigma_\rho]=n^{\downarrow |\rho|}\,(1-q)^{|\rho|}/(1-q^\rho)$.
\bigskip

\subsection{Bound for cumulants}\label{subsect:bound_cumulants} This sort of result (computing the fluctuations of an observable from the character values) has already been stated by P. \'Sniady in the case where the characters are small and the corresponding diagrams balanced. Here, we have to use different observables and a different gradation, but the ideas of the proof are the same, and the definitions of the different cumulants used here come from \cite{Sniady06}. In particular we shall use joint cumulants of random variables $k(X_1,\ldots,X_r)$, which can be defined by induction\footnote{It can also be defined via generating functions.} on $r$:

$$\esper[X_1 \,X_2 \cdots X_r]=\sum_{\substack{\pi \text{ partition}\\ \text{of } \{1,2,\ldots,r\} }} \, k\big(X_{i \in \pi_1}\big) \, k\big(X_{i \in \pi_2}\big) \cdots k\big(X_{i \in \pi_\ell}\big). $$

\noindent A centered gaussian vector $(Y_1,\ldots,Y_k)$ is characterized by the fact that all $k(Y_{i_1},\ldots,Y_{i_r})$ are equal to $0$ for $r \neq 2$. Moreover, the values for $r=2$ give then the covariance matrix. Therefore, the proof of Proposition \ref{prop:deviation_powersums} consists in computing the limit of $k(W_{n,q,i_1},\ldots,W_{n,q,i_r})$.
In fact, we will show a stronger result than the convergence towards $0$ for $r \geq 3$:
\begin{lemma}[Order of magnitude of cumulants]\label{lem:ordre_cumulants}
For any $x_1,x_2,\ldots,x_r \in \obs$, one has:
$$k(x_1,\ldots,x_r) = O\left(n^{\deg(x_1) + \cdots + \deg(x_r)-r+1}\right).$$
\end{lemma}\bigskip

The end of this subsection is devoted to the proof of this lemma.
Thanks to the multi-linearity of cumulants (for $r\geq 2$) and some compatibility relation with respect to products (see \cite{LS59} and \cite[Theorem 4.4]{Sniady06}), it is enough to prove this lemma with the $x_i$ in any graded algebraic basis of $\obs$. A good choice is the basis $(\varSigma_k)_{k\geq 1}$; indeed, we have already computed the expectation of these observables under the $q$-Plancherel measure. Unfortunately, we do not know their joint moments $\esper[\varSigma_{i_1} \cdots \varSigma_{i_r}]$, which are different from what we have computed in the previous section, namely, $\esper[\varSigma_{i_1,\ldots,i_r}]$.\bigskip

\subsubsection{Cumulants $k^\bullet$}
This difficulty disappears if we use the algebra $\obs^\bullet$, because $\Sigma_{i_1} \bullet \cdots \bullet \Sigma_{i_r}=\Sigma_{i_1,\ldots,i_r}$. This algebra gives rise to a new kind of cumulant, the so-called disjoint cumulants, which we shall denote by $k^\bullet(X_{1},\ldots,X_{r})$. Their definition is the same as for $k(X_1,\ldots,X_r)$, except that the product of observables is now the disjoint product:
\begin{equation}\label{eq:def_dis_cum}
\esper[X_1 \bullet X_2 \bullet \cdots \bullet X_r]=\sum_{\substack{\pi \text{ partition}\\ \text{of } \{1,2,\ldots,r\} }} \, k^\bullet\big(X_{i \in \pi_1}\big) \, k^\bullet\big(X_{i \in \pi_2}\big) \cdots k^\bullet\big(X_{i \in \pi_\ell}\big).
\end{equation}
It is easy to prove the analog of Lemma \ref{lem:ordre_cumulants} for these disjoint cumulants:
\begin{lemma}[Order of magnitude of disjoint cumulants]\label{lem:ordre_disjoint_cumulants}
For any $x_1,x_2,\ldots,x_r \in \obs$, one has:
$$k^\bullet(x_1,\ldots,x_r) = O\left(n^{\deg(x_1) + \cdots + \deg(x_r)-r+1}\right).$$
\end{lemma}\bigskip
\begin{proof}
As for classical cumulants, it is enough to prove it for $x_j=\varSigma_{i_j}$. In this case, the left-hand side of \eqref{eq:def_dis_cum} is simply
$$\esper[\varSigma_{i_1,\ldots,i_r}]=\frac{(q-1)^{i_1+\cdots+i_r}}{\prod(q^{i_j}-1)} n^{\downarrow (i_1+\cdots+i_r)}.$$
Except for the factor $n^{\downarrow (i_1+\ldots+i_r)}$, everything is multiplicative with respect to the parts of $\lambda$. More formally, there exist numbers $\alpha_i$ such that:
\begin{align*}\esper[\varSigma_{i_1,\ldots,i_r}] &= \left( \prod_j \alpha_{i_j} \right) n^{\downarrow (i_1+\cdots+i_r)};\\
&=\left( \prod_j \alpha_{i_j} \right) \esper\big[\varSigma_{\!\!\!\!\underbrace{1,\ldots,1}_{i_1 + \cdots + i_r \text{ times}}}\big].
\end{align*}
An immediate induction implies that:
$$k^\bullet(\varSigma_{i_1},\ldots,\varSigma_{i_n})=\left( \prod_j \alpha_{i_j} \right) k^\bullet(\varSigma_{\underbrace{1,\ldots,1}_{i_1 \text{ times}}},\ldots,\varSigma_{\underbrace{1,\ldots,1}_{i_r \text{ times}}}).$$
But $\chi^\lambda(\varSigma_{1^{i}})=n^{\downarrow i}$ for any $\lambda$, so the corresponding cumulants do not depend on the measure considered on Young diagrams, and we can use the result of P. \'Sniady who proved (\cite[Lemma 4.8]{Sniady06}) that:
$$k^\bullet(\varSigma_{\underbrace{1,\ldots,1}_{i_1 \text{ times}}},\ldots,\varSigma_{\underbrace{1,\ldots,1}_{i_r \text{ times}}})=O(n^{i_1+\cdots+i_r-(r-1)}).$$
This ends the proof of this lemma.
\end{proof}

\subsubsection{Cumulants $k^{\id}$}
The last problem is to link classical and disjoint cumulants. Thanks to the commutativity of the diagram below
$$\begin{diagram}
   \node{\obs} \arrow{se,b}{\esper=M_{n,q}}\arrow[2]{e,t}{\esper^\id} \node[2]{\obs^\bullet}\arrow{sw,b}{\esper=M_{n,q}}\\
\node[2]{\C}
  \end{diagram}
$$
this can be done by using a result of Brillinger \cite{Bri69}. One obtains the following formula, see \cite[Proposition 4.1]{Sniady06}:
\begin{equation}\label{eq:link_2_cumulants}
k(X_1,X_2,\ldots,X_n)=\sum_{\pi\vdash n } k^\bullet(k^{\id}(X_{i \in \pi_1}), k^{\id}(X_{i \in \pi_2}),\ldots,k^{\id}(X_{i \in \pi_r})),
\end{equation}
where $k^{\id}$ is the cumulant corresponding to the identity map $\esper^\id$ of the diagram (note that $k^{\id}$ takes value in $\obs$), that is to say that it is defined by
$$
X_1 X_2 \cdots X_r=\sum_{\substack{\pi \text{ partition}\\ \text{of } \{1,2,\ldots,r\} }} k^{\id}\big(X_{i \in \pi_1}\big) \bullet k^{\id}\big(X_{i \in \pi_2}\big) \bullet \cdots \bullet k^{\id}\big(X_{i \in \pi_{\ell}}\big).
$$

We will prove the following bound on the degree of identity cumulants. For any $i_1,\ldots,i_r$, one has:
\begin{equation}\label{eq:major_deg_kid}
\deg\big(k^{\id}(\Sigma_{i_1},\ldots,\Sigma_{i_r})\big) \leq i_1 + \cdots + i_r - (r-1).
\end{equation}
Note that it is not the same as \cite[Theorem 4.3]{Sniady06}; indeed, we use a different gradation on the algebra of observables.\bigskip

Let us begin by looking at $k^{\id}(\varSigma_{i_1},\ldots,\varSigma_{i_n})$ for small values of $n$.
For $n=1$, by definition $k^{\id}(\varSigma_{i})=\varSigma_i$. For $n=2$, one has:
$$k^{\id}(\varSigma_{i_1},\varSigma_{i_2})=\varSigma_{i_1} \cdot \varSigma_{i_2} - \varSigma_{i_1,i_2}.$$
But one also has
\begin{align*}
 \varSigma_{i_1} \cdot \varSigma_{i_2} &=\sum_{ A_1\in \Arr(n,i_1),\,A_2 \in \Arr(n,i_2)} C(A_1) \, C(A_2) ;\\
\varSigma_{i_1,i_2} &= \sum_{\substack{ A_1\in \Arr(n,i_1),\,A_2 \in \Arr(n,i_2)\\ A_1 \cap A_2 = \emptyset}} C(A_1)\, C(A_2).
\end{align*}
so the difference $k^{\id}$ is equal to the same expression with the sum restricted to pairs of arrangements which have a non-empty intersection:
$$k^{\id}(\varSigma_{i_1},\varSigma_{i_2} )=\sum_{ \substack{A_1\in\Arr(n,i_1),\,A_2 \in \Arr(n,i_2)\\ A_1 \cap A_2 \neq \emptyset}} C(A_1)\, C(A_2).$$
It turns out that the identity cumulants satisfy a similar identity in the general case, and the intersections between the $A_j$ play once again a significant role. First, note that:
$$\varSigma_{i_1} \cdots  \varSigma_{i_r} =\sum_{ \forall j,\,\,A_j \in \Arr(n,i_j)} C(A_1) \cdots C(A_r)$$
Given a family $A_1,\ldots,A_r$ of arrangements, we consider the relation $\sim_1$ defined by $k \sim_1 l \Leftrightarrow A_k \cap A_l \neq \emptyset$. Its transitive closure is an equivalence relation; let us denote by $\pi(A_1,\ldots,A_r)$ the corresponding partition of $\lle 1,r\rre$. Obviously, this partition depends only on which indices are equal and not on the values they take, so the sum of $C(A_1) \cdots C(A_r)$ over arrangements corresponding to a given partition is of the form $X_{i_1,\dots,i_r}^{\mathcal{G}}$ (notation of paragraph \ref{subsect:bullet_product}) and we can use Lemma \ref{lem:bullet} on this kind of expression.

\begin{lemma}
For any $i_1,\ldots,i_r$, one has:
$$k^{\id}(\varSigma_{i_1},\ldots,\varSigma_{i_r}) = \sum_{ \substack{\forall j,\,\,A_j \in \Arr(n,i_j) \\  \pi(A_1,\ldots,A_r)=\{\lle 1,r\rre\}       }} C(A_1) \cdots C(A_r).$$
\end{lemma}

\begin{proof}
We will prove this statement by induction on $r$. As explained above, it is true for $r=1,2$.
Let us suppose that it is true for all $s \leq r-1$. Because of the definition of the disjoint product $\bullet$, for every \emph{non-trivial} set partition $\pi$, one has: 
$$k^{\id}\big(\varSigma_{i \in \pi_1}\big) \bullet k^{\id}\big(\varSigma_{i \in \pi_2}\big) \bullet \cdots \bullet k^{\id}\big(\varSigma_{i \in \pi_\ell}\big)=\sum_{\substack{\forall j,\,\,A_j \in \Arr(n,i_j) \\  \pi(A_1,\ldots,A_r)=\pi} } C(A_1) \cdots C(A_r).$$
Indeed, the induction hypothesis applies for each part $\pi_i$, so we only keep terms corresponding to sequences of arrangements for which two numbers in the same part of $\pi$ are in the same part of $\pi(A_1,\ldots,A_r)$.
Besides, when we form the product $\bullet$, we take only disjoint set of indices (using Lemma \ref{lem:bullet}), so two numbers in different parts of $\pi$ cannot be in the same part of $\pi(A_1,\ldots,A_r)$.\bigskip

\noindent We use now the inductive definition of $k^{\id}$:
\begin{align*}
k^{\id}(\varSigma_{i_1},\ldots,\varSigma_{i_r}) &= \,\,\,\,\varSigma_{i_1} \cdots \varSigma_{i_r}\quad  - \sum_{\pi \neq \{\lle1,r\rre\}} \, k^{\id}\big(\varSigma_{i \in \pi_1}\big) \bullet k^{\id}\big(\varSigma_{i \in \pi_2}\big) \bullet \cdots \bullet k^{\id}\big(\varSigma_{i \in \pi_\ell}\big) \\
&=\sum_{ \forall j,\,\,A_j\in \Arr(n,i_j)} \!\!\!\!C(A_1) \cdots C(A_r) \quad - \sum_{\pi \neq \{\lle 1,r\rre \}} \left( \sum_{ \substack{ \forall j,\,\,A_j \in \Arr(n,i_j)\\ \pi(A_1,\ldots,A_r) =\pi}} C(A_1) \cdots C(A_r) \right)\\
& =\!\!\!\!\sum_{ \substack{\forall j,\,\,A_j \in \Arr(n,i_j) \\  \pi(A_1,\ldots,A_r)=\{\lle 1,r\rre\}}} \!\!\!\!\!\!\!\!C(A_1) \cdots C(A_r),\end{align*}
which is exactly the equality we wanted at rank $r$.
\end{proof}
Now, for every sequence of arrangements $A_1,\ldots,A_r$ such that $\pi(A_1,\ldots,A_r)$ is the trivial partition $\{\lle 1,r \rre\}$, the cardinality of $|A_1 \cup \cdots \cup A_r|$ is smaller than $i_1 + \cdots + i_r - (r-1)$, because one has at least $r-1$ equalities between different elements of these arrangements. So, $k^{\id}(\varSigma_{i_1},\ldots,\varSigma_{i_r})$ can be written as a linear combination of $\varSigma_{\lambda}$ with $|\lambda| \leq i_1 + \cdots + i_r - (r-1)$, which ends the proof of equation \eqref{eq:major_deg_kid}.\bigskip

Finally, equations \eqref{eq:link_2_cumulants} and \eqref{eq:major_deg_kid}, together with Lemma \ref{lem:ordre_disjoint_cumulants}, complete the proof of Lemma \ref{lem:ordre_cumulants}.
\bigskip

\subsection{Computation of covariances} Now, let us compute explicitly the limit (co)variances of the rescaled deviations of the $\varSigma_k$. First, note that if $l$ and $m$ are two positive integers, then:
\begin{align*}
k^\bullet(\varSigma_l,\varSigma_m)&=\esper[\varSigma_l \bullet \varSigma_m] - \esper[\varSigma_l]\,\esper[\varSigma_m] = \esper[\varSigma_{l,m}] - \esper[\varSigma_l]\,\esper[\varSigma_m]\\
&= \frac{(1-q)^{l+m}}{1-q^{l,m}}\,\left(n^{\downarrow l+m} - n^{\downarrow l}\,n^{\downarrow m}\right)=-lm\,\frac{(1-q)^{l+m}}{1-q^{l,m}}\,n^{l+m-1}+ O(n^{l+m-2})\end{align*}
The last identity follows from the development of the falling factorial:
$$n^{\downarrow k} =n^k- \frac{k(k-1)}{2}\,n^{k-1}+O(n^{k-2})$$
Now, we can compute the higher term of the standard cumulant $k(\varSigma_l,\varSigma_m)$ if we use the graded development of $\varSigma_l\,\varSigma_m-\varSigma_{l,m}$ in the algebra of observables:
\begin{align*}k(\varSigma_l,\varSigma_m)&=k^{\mathrm{id}}(\varSigma_l,\varSigma_m)+k^\bullet(\varSigma_l,\varSigma_m)=\esper[\varSigma_{l}\,\varSigma_{m}-\varSigma_{l,m}]+k^\bullet(\varSigma_l,\varSigma_m)\\
&=\esper[lm\,\varSigma_{l+m-1}]-lm\,\frac{(1-q)^{l+m}}{1-q^{l,m}}\,n^{l+m-1}+ O(n^{l+m-2})\\
&=lm\,(1-q)^{l+m}\,\left(\frac{1}{1-q^{l+m-1,1}} - \frac{1}{1-q^{l,m}}\right)\,n^{l+m-1}+ O(n^{l+m-2})
\end{align*}
We have therefore proved:
\begin{proposition}[Fluctuations of $\varSigma_k$]\label{PropFluctuationsCharacters}
If we denote $Z_{n,q,k}=\sqrt{n}\,(\varSigma_k/n^k - (1-q)^k/(1-q^k))$, then $(Z_{n,q,k})_{k \geq 2}$ converges in finite-dimensional laws towards a gaussian vector $(Z_{q,k})_{k \geq 2}$, with the following covariance matrix:
$$k(Z_{q,l},Z_{q,m})=lm\,(1-q)^{l+m}\,\left(\frac{1}{1-q^{l+m-1,1}} - \frac{1}{1-q^{l,m}}\right)$$
\end{proposition}
\begin{proof}[End of the proof of Proposition  \ref{prop:deviation_powersums}] It turns out that the sequence $(W_{n,q,k})_{k\geq 2}$ of the rescaled deviations of power sums has the same limit law $(Z_{q,k})_{k \geq 2}$. Indeed, we have already shown that the limit of the rescaled deviations of power sums is a gaussian vector, and it suffices then to show that $k(p_l,p_m)$ has the same higher term as $k(\varSigma_l,\varSigma_m)$. According to the proposition \ref{prop:sigma_pk}, each symbol $\varSigma_k$ is a rational polynomial in the $p_{l\leq k}$ with higher term $p_k$.  So, we can use the multi-linearity of joint cumulants to write:
$$k(\varSigma_l,\varSigma_m)=k(p_l,p_m) + (\text{cumulants of observables }a_i,b_j\text{ with }\deg a_i+\deg b_j <l+m).$$
But if $a$ and $b$ are two observables of degree $A$ and $B$, then we have seen that $k(a,b)$ is a $O(n^{A+B-1})$. As a consequence, the remaining terms are all $O(n^{l+m-2})$. This ends the proof of Proposition \ref{prop:deviation_powersums}.
\end{proof}\medskip

\begin{remark}\label{RemFluctuationsQCharacters}
Lemma \ref{lem:ordre_cumulants} also implies that the rescaled $q$-characters converge towards a gaussian vector. For these $q$-characters, it can be shown that if $S_{n,q,k}=\sqrt{n}\,(\varSigma_{k,q}/n^{k})\simeq \sqrt{n} \,\chi^{\lambda}(k1^{n-k},q)$, then $(S_{n,q,k})_{k \geq 2}$ converges in finite-dimensional laws towards a gaussian vector $(S_{q,k})_{k \geq 2}$, with the following covariance matrix:
$$k(S_{q,l},S_{q,m})=(q-q^{2})^{l+m-3}\,(1-q^{2})\,\frac{\{l-1\}_{q}\,\{m-1\}_{q}}{\{l+m-1\}_{q}\,\{l+m-2\}_{q}\,\{l+m-3\}_{q}}$$
Hence, one has a $q$-analogue of Kerov's central limit theorem (see \cite[Theorem 6.1]{IO02}) for the $q$-characters of the Hecke algebras under $q$-Plancherel measures. However, notice that the fluctuations of the $q$-characters are always some $O(n^{-1/2})$, in contrast to the case of the characters $\chi^{\lambda}(k1^{n-k})$ under Plancherel measures, whose fluctuations are $O(n^{-k/2})$. The proof of the formula given above involves a generalization of Proposition \ref{PropFluctuationsCharacters} to the case of arbitrary central characters $\varSigma_{\mu}$, and a M\"obius  inversion formula that allows to simplify the expression of the covariances; we refer to a future paper of the second author (\cite{Mel10}) for details on these computations.
\end{remark}\bigskip

\subsection{Fluctuations of rows} \label{subsect:fluctuations_row}
 As before, we denote by $X_\lambda$ the probability measure on $[-1,1]$ associated to a Young diagram $\lambda$, and by $X_{\infty;q}$ the limit probability measure of the $X_{\lambda}$ under the $q$-Plancherel measures. If $f \in \mathscr{C}^1([-1,1])$ is a continuously differentiable function on $[-1;1]$, we define the rescaled deviation of $f$ under $q$-Plancherel measure. It is the random variable $D_{n,q}(f)$ obtained by composition of $q$-Plancherel measure by the function:
$$\lambda \mapsto \sqrt{n}\,\left(X_\lambda(f)-X_{\infty;q}(f)\right),$$
where $\mu(f)$ is a short notation for $\int_{-1}^1 f d\mu$.\bigskip

In this context, the convergence in law of rescaled power sums can be reformulated as a convergence in law of the vector 
$\big(D_{n,q}(x^l)\big)_{1 \leq l \leq k}$ towards a gaussian vector $(Z_{q;l})$ with covariance matrix: 
\begin{align*}k(Z_{q;l},Z_{q;m}) &= X_{\infty;q}((l+1)(m+1)\,x^{l}x^m)-(l+1)(m+1)\,X_{\infty;q}(x^{l})\,X_{\infty;q}(x^{m}) ;\\
&=X_{\infty;q}\left((x^{l+1})'\,(x^{m+1})'\right)-X_{\infty;q}\left((x^{l+1})'\right)\,X_{\infty;q}\left((x^{m+1})'\right).
\end{align*}
Recall that with respect to the $\mathscr{C}^1$-topology, $\R[X]$ is dense in the space of continuously differentiable functions $\mathscr{C}^1([-1,1])$.
By using this density together with the bilinearity of $D_{n,q}$ and covariance, we conclude that for any functions $f_1,f_2,\ldots,f_r \in \mathscr{C}^1([-1,1])$, the vector of rescaled deviations
$D_{n,q}(f_i)$ converges towards a gaussian vector $(Z_{q;f_i})_{i \in \lle 1,r\rre}$ with covariance matrix:
$$k(Z_{q;f_i},Z_{q;f_j}) = X_{\infty;q}\left((xf_i)'\,(xf_j)'\right)-X_{\infty;q}\left((xf_i)'\right)\,X_{\infty;q}\left((xf_j)'\right).$$

\begin{proof}[Proof of Theorem \ref{th:second_order_as_rows}] If $i \geq 1$, we denote by $f_i$ a class $\mathscr{C}^1$ and non-negative function such that $f_i(x)=1$ on a vicinity $V_i$ of $q^{i-1}\,(1-q)$ and $f_i(x)=0$ outside a vicinity $W_i\supset V_i$, \emph{cf.} figure \ref{galactica}. We shall also suppose that the $W_i$ are disjoint open sets; in particular, $X_{\infty;q}(f_i)=q^{i-1}(1-q)$.  Since $f_i'(q^{i-1}(1-q))=0$, $(xf_i)'(q^{i-1}(1-q)) = f_i(q^{i-1}(1-q))=1$, and as a consequence,
$$X_{\infty;q}((x\,f_i(x))'^2)-X_{\infty;q}((x\,f_i(x))')^2=q^{i-1}(1-q)-q^{2(i-1)}(1-q)^2.$$

\figcap{\psset{unit=1mm}\pspicture(-20,-4)(80,13)
\psline{->}(-20,0)(80,0)
\psline{->}(0,-2)(0,20)
\pscurve[linecolor=blue](-15,0.2)(5,0.2)(20,0.2)(34,0.2)(39,1)(46,14)(50,15)(54,14)(61,1)(67,0.2)(75,0.2)
\psline(75,-1)(75,1)
\psline(50,-1)(50,1)
\psline(37,-8)(37,1)
\psline(63,-8)(63,1)
\psline{<->}(37,-7)(63,-7)
\rput(50,-9.5){$W_i$}
\rput(50,3){$q^{i-1}(1-q)$}
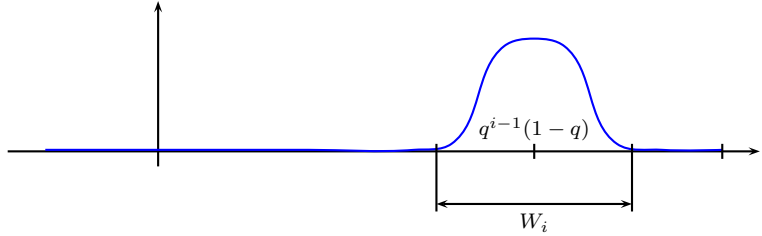
\endpspicture}{Graphical representation of the test function $f_i$.\label{galactica}}

Moreover, $f_i\,f_j=0$ for all $i\neq j$, and the same holds for $f_i\,f_j'$, $f_i'\,f_j$ or $f_i'\,f_j'$. Therefore, the covariance matrix of the corresponding limit gaussian vector $(Z_{q;f_i})_{i \geq 1}$ is:
\begin{align*}X_{\infty;q}\left((x\,f_i(x))'\,(x\,f_j(x))'\right)-X_{\infty;q}\left((x\,f_i(x))'\right)\,X_{\infty;q}\left((x\,f_j(x))'\right)&=-X_{\infty;q}\left((x\,f_i(x))'\right)\,X_{\infty;q}\left((x\,f_j(x))'\right) ;\\
&=-(1-q)^2\,q^{i+j-2}.\end{align*}
Finally, for $n$ large enough and outside a set of arbitrary small probability, $a_i^*(\lambda)/n$ is in $V_i$, $a_{i-1}^*(\lambda)/n$ is in $V_{i-1}$ and $a_{i+1}^*(\lambda)/n$ is in $V_{i+1}$. Thus, with a probability close to $1$, the support of $X_\lambda$ and the support of $f_i \neq 0$ intersect only at $a_i^*(\lambda)/n$, and $X_\lambda(f_i)=a_i^*(\lambda)/n$. As a consequence, $$\sqrt{n}\,\left(\frac{a_i^*(\lambda)}{n}-q^{i-1}(1-q)\right)=\sqrt{n}\,\left(X_\lambda(f_i)-X_{\infty;q}(f_i)\right),$$  
outside a set of arbitrary small probability, and the repartition functions of these random variables have the same asymptotic behavior. Theorem \ref{th:second_order_as_rows} follows then from the estimation $a_i^*(\lambda)=\lambda_i+O(1)$.\end{proof}
\bigskip

\section{Schur-Weyl representations}\label{sect:schur-weyl}
In this section, we use the same tools to study other measures on partitions and to describe asymptotically the shape of the diagrams. The purpose of this section is to show that our method works as soon as the typical length of the first rows of the diagram is of order $n^\alpha$ with $\alpha>1/2$ and not only with $\alpha=1$.

\subsection{Decomposition of tensor representations}\label{subsect:def_SW}
We consider the following representation of the symmetric group: $\sym_n$ acts on the space $(\C^N)^{\otimes n}$ by commuting the vectors of the pure tensors. Its character is very easy to compute: $\chi(\rho)=N^{\ell(\rho)}$, where $\ell(\rho)$ is the number of parts of the partition $\rho$ (or, equivalently, the number of cycles of any permutation of type $\rho$).\bigskip

This representation can be decomposed into irreducible representations of $\sym_n$: let us denote $m_\lambda$ the multiplicity of the irreducible module $S_\lambda$. We associate to this decomposition the following measure on partitions
$$\proba[\lambda] = \frac{m_\lambda \cdot \dim(S_\lambda)}{N^n}.$$
In other words, $\proba[\lambda]$ is the dimension of the isotypic component of type $\lambda$ divided by the total dimension of the module $(\C^N)^{\otimes n}$. Note that this is a very natural way to construct a measure on Young diagrams of a given size: the same construction with the regular representation gives  rise to the Plancherel measure.\bigskip

What is interesting with this kind of measure is that the expectation of rescaled character values is proportional to the character of the representation. In the case we are looking at (the action of $\sym_n$ on $(\C^N)^{\otimes n}$),
\begin{align*}
 \esper[\Sigma_\rho] & = n^{\downarrow |\rho|} \cdot \frac{\chi(\rho 1^{n-|\rho|})}{N^n} \\ &= n^{\downarrow |\rho|} \cdot N^{\ell(\rho) - |\rho|}.
\end{align*}
Let us finish this paragraph with a little remark. By Schur-Weyl duality, one has the following isomorphism of $\sym_n \times \GL(N,\C)$-modules
$$(\C^N)^{\otimes n} \simeq \bigoplus_{\lambda \vdash n,\ell(\lambda)\leq N} V_\lambda \otimes U_\lambda,$$
where $U_\lambda$ is an irreducible $\GL(N,\C)$-module. Moreover, the dimension of $U_\lambda$ is the number of semi-standard tableaux of shape $\lambda$ with entries lower or equal than $N$. So,
$$\proba[\lambda] = \frac{\left|\left\{\text{\begin{tabular}{c}
                         standard tableaux\\
			 of shape $\lambda$
                        \end{tabular}}\right\}\right|
				\cdot
\left|\left\{\text{\begin{tabular}{c}
                         semistandard tableaux\\
			 of shape $\lambda$ (entries $\leq N$)
                        \end{tabular}}\right\}\right|
}{N^n}.$$
This implies, as mentioned in paragraph \ref{subsect:ss_suite_croissante}, that this measure is in fact the image of the uniform distribution on words of length $n$ and letters from $1$ to $N$ by the RSK algorithm. On the other hand, notice that the Schur-Weyl measure charges only the Young diagrams whose lengths are smaller than $N$. \bigskip

\subsection{Asymptotic shape of the diagrams}
In the previous paragraph, we introduced some measures on Young diagrams of size $n$. Here we are interested in the asymptotic shape of typical diagrams when $n$ tends to infinity. That means, as the measure we introduced also depends on an other parameter $N$, that we have to choose a value of $N$ for each $n$. The case where the parameter $N$ is equivalent of an expression of the kind $1/c \cdot n^\alpha$ with $\alpha \geq 1/2$ has already been solved by P. Biane \cite[paragraph 3.1]{Bia01}. So we will look at the case $\alpha < 1/2$.\bigskip

Using the formula for expectation of normalized character values of the previous paragraph, one has:
\begin{equation}\label{eq:equiv_Esigma}
\esper[\varSigma_\rho] = n^{\downarrow |\rho|} \cdot N^{\ell(\rho)-|\rho|} \sim c^{|\rho| - \ell(\rho)} n^{|\rho| + \alpha \ell(\rho) -\alpha |\rho|}.
\end{equation}
As in section \ref{sect:first_order}, the first step is to show that the power sums have the same asymptotic behavior. Lemma \ref{lem:major_Ex_deg} is true in this context, but not strong enough to obtain such a result. Unfortunately, with the gradation we give in section \ref{sect:observables}, it is not possible to have a stronger upper bound. Thus one has to introduce another gradation on $\obs$.\bigskip

\subsubsection{Gradations on $\obs$}
Let us define (this may not be an integer, but it does not matter):
$$\deg_\alpha(p_\rho) = |\rho| + \alpha \ell(\rho) -\alpha |\rho|.$$
As in the case $\alpha=0$, the expression of $\varSigma_k$ in terms of power sums (Proposition \ref{prop:sigma_pk}) gives us the degree and the top homogeneous component of $\varSigma_\rho$.
\begin{lemma}\label{lem:prho_encore_dominant}
 For any $\rho \in \young$, one has:
$$\varSigma_\rho = p_\rho + \text{terms of lower degree}.$$
\end{lemma}
\begin{proof}
It is enough to prove this for $\rho=(k)$. Note that the $\alpha$-degree can be written as:
$$\deg_\alpha(p_\lambda)= (1 - 2 \alpha) |\lambda| + \alpha (|\lambda| + \ell(\lambda)).$$
\begin{enumerate}
\item The first term is a positive multiple (recall that $\alpha<1/2$) of the usual gradation, so if $p_\lambda$ appears with a non-zero coefficient in $\varSigma_k$, one has $|\lambda| \leq k$, with equality if and only if $\lambda=(k)$.
\item The second term is a non-negative multiple of $2 \deg_{1/2}(\lambda)=|\lambda| + \ell(\lambda)$. An easy consequence of Proposition \ref{prop:sigma_pk} is that all $p_\lambda$ appearing in $\varSigma_k$ fulfill the condition $|\lambda| + \ell(\lambda) \leq k+1$. Note that we do not have the same necessary condition of equality as before: in fact, $\deg_{1/2}$ is a classical gradation of $\obs$ (see \cite{IO02}), for which the top homogeneous component of $\varSigma_\rho$ \emph{is not} $p_\rho$.
\end{enumerate}
Finally, $\varSigma_\rho$ has degree $|\rho| + \alpha \ell(\rho) -\alpha |\rho|$, and its top homogeneous component for the new gradation is $p_\rho$.
\end{proof}\bigskip

\subsubsection{Convergence of rescaled power sums and concentration of characters}
Using equation (\ref{eq:equiv_Esigma}), one has the following analog of Lemma \ref{lem:major_Ex_deg} (the gradation has been chosen for this purpose):
\begin{lemma}\label{lem:Ex_encore_majore}
  For any $x \in \obs$, one has: 
$\esper[x] = O(n^{\deg_\alpha(x)})$.
\end{lemma}

\noindent Now, we have all the tools to prove the convergence of rescaled powers sums:

\begin{proposition}\label{prop:cv_powersums_encore}
Let $(N_n)_{n \geq 1}$ be a sequence of positive integers such that $N_n \sim c \cdot n^\alpha$ for some $\alpha <1/2$. Under the measure associated with the action of the symmetric groups on the spaces $(\C^{N_n})^{\otimes n}$, one has the following convergence in probability:
$$\frac{p_k}{n^{k + \alpha - \alpha k}} \rightarrow c^{k-1}.$$
\end{proposition}

\begin{proof}
 Equation (\ref{eq:equiv_Esigma}), Lemma \ref{lem:prho_encore_dominant} and Lemma \ref{lem:Ex_encore_majore} imply that:
$$\esper[p_\rho] \sim c^{\rho - \ell(\rho)} n^{|\rho| + \alpha \ell(\rho) -\alpha |\rho|}.$$
The result follows thanks to Bienaym\'e-Chebyshev inequality as in the case $\alpha=0$.
\end{proof}

\begin{corollary}[Concentration of characters]
 Under the same assumptions, in probability,
$$\frac{\varSigma_\rho}{n^{|\rho| + \alpha \ell(\rho) -\alpha |\rho|}} \rightarrow c^{|\rho|-\ell(\rho)}.$$
\end{corollary}\medskip

\subsubsection{Consequences on the shape of diagrams}
To deduce some results on the shape of diagrams from the convergence of the quantities $p_k/(n^{k + \alpha - \alpha k})$, one has to interpret them as the moments of some measure-valued random variables. Therefore we define for a Young diagram $\lambda$ the measure:
$$X_{\lambda,\alpha}=\sum_{i=1}^d (a_i^*(\lambda)/n^{1-\alpha}) \, \frac{\delta_{(a_i^*(\lambda)/n^{1-\alpha})}}{n^\alpha} + (b_i^*(\lambda)/n^{1-\alpha}) \, \frac{\delta_{(-b_i^*(\lambda)/n^{1-\alpha})}}{n^\alpha}.$$
With this definition, Proposition \ref{prop:cv_powersums_encore} can be read as the convergence in probability of the moments of this random-valued measure towards the sequence $(1,c,c^2,\ldots)$, which is obviously the sequence of the moments of the measure $\delta_c$ (the Dirac mass at $c$). Since the probability measure $\delta_{c}$ is characterized by its moments, this implies the convergence in law $X_{\lambda,\alpha} \to \delta_{c}$ in probability in the space of probability measures $\mathscr{M}(\R)$. \bigskip

Now, the negative part of $X_{\lambda,\alpha}$ is in fact negligible. Indeed, since $\ell(\lambda)\leq N$ with probability $1$, the total weight of 
$$\sum_{i=1}^{d} (b_{i}^{*}(\lambda)/n^{1-\alpha})\,\frac{\delta_{(-b_{i}^{*}(\lambda)/n^{1-\alpha})}}{n^{\alpha}}$$
is bounded from above and with probability $1$ by:
$$\sum_{i=1}^{N}\frac{N}{n} = \frac{N^{2}}{n}\simeq c^{-2}\,n^{2\alpha-1}\to 0$$
As a consequence, one can replace $X_{\lambda,\alpha}$ by its positive part $\sum_{i=1}^{d}(a_{i}^{*}(\lambda)/n^{1-\alpha})\,\frac{\delta_{(a_{i}^{*}(\lambda)/n^{1-\alpha})}}{n^{\alpha}}$, and even by the simpler probability measure 
$$Y_{\lambda,\alpha}=\sum_{i=1}^{N} (\lambda_{i}/n^{1-\alpha})\,\frac{\delta_{\lambda_{i}/n^{1-\alpha}}}{n^{\alpha}}=\sum_{i=1}^{N} \frac{\lambda_{i}}{n}\delta_{\lambda_{i}/n^{1-\alpha}};$$
this random measure converges in probability towards $\delta_{c}$ in $\mathscr{M}(\R)$.\bigskip

We fix some positive real numbers $\theta$, $\eta$ and $\eps$, and we consider a continuous function $f:\R \to [0;1]$ such that $f(x)=0$ outside the interval $[c(1-\eta/2),c(1+\eta/2)]$, and $f(c)=1$. Then, for $n$ big enough and with probability greater than $1-\eps$,
$$1-\theta=f(c)-\theta \leq \sum_{i=1}^{N} \frac{\lambda_{i}}{n}\,f\left(\frac{\lambda_{i}}{n^{1-\alpha}} \right) \leq f(c)+\theta=1+\theta$$
because of the convergence in law $Y_{\lambda,\alpha}\to \delta_{c}$. The sum involved in this estimate is bounded from above by $\sum_{\left\{i\,\,|\,\,\lambda_{i}/n^{1-\alpha} \in [c(1-\eta/2),c(1+\eta/2)]\right\}}\frac{\lambda_{i}}{n}$, and therefore by
$$\frac{c(1+\eta/2)}{n^{\alpha}}\,\,\card\left\{i\,\,\bigg|\,\,\frac{\lambda_{i}}{n^{1-\alpha}} \in [c(1-\eta/2),c(1+\eta/2)]\right\}.$$
As a consequence, the number of rows $\lambda_{i}$ such that $c(1-\eta/2)\leq \lambda_{i}/n^{1-\eta} \leq c(1+\eta/2)$ is bigger than 
$$\frac{1-\theta}{c(1+\eta/2)} \,n^{\alpha}\geq \frac{(1-\theta)(1-\eta/2)}{c}\,n^{\alpha}\geq \frac{1-\eta}{c}\,n^{\alpha}$$
for $\theta$ small enough.
As by the very definition of our measure, our random diagrams can not have
columns of length bigger than $N$, the proof of Theorem~\ref{prop:SW_asymptotic_shape} is over.

\begin{remark}
The vector space $(\C^N)^{\otimes n}$ is also a representation of the unitary group $U_N$ and, using Schur-Weyl duality, its isotypic components under the action of $U_N$ are the same as under the action of $S_n$. As a consequence, the results of article \cite{CollinsSniady09} can be used to study the asymptotic shape of the diagram. It turns out that, in the case $\alpha <1/2$, this approach gives finer results than ours. Nevertheless, as they are two independent parameters, the space $(\C^N)^{\otimes n}$ remains a good example to show that our method works with any intermediate decay of character values between $O(1)$ and $O(\sqrt{n}^{-|\rho|+\ell(\rho)})$.
\end{remark}

\section*{Acknowledgements}
The authors would like to thank P. Biane for his valuable help, as well as A. Borodin and P. \'Sniady for interesting remarks. In addition, we want to mention that Theorem \ref{th:first_order_as_rows} was conjectured by the second author from examples (like Figure \ref{fig:ex_diag_qplanch}) obtained by computer exploration using the open-source
mathematical software \texttt{Sage}~\cite{sage} and its algebraic combinatorics features developed by the \texttt{Sage-Combinat}
community~\cite{Sage-Combinat}.\bigskip

Finally, we are indebted to anonymous referees, whose comments have been very useful to make our development clearer and more comprehensive.

\bibliographystyle{alpha}
\bibliography{qplancherel}

\newcommand{\etalchar}[1]{$^{#1}$}
\begin{thebibliography}{RRW96}

\bibitem[BDJ99]{BDJ99firstrow}
J.~Baik, P.~Deift, and K.~Johansson.
\newblock On the distribution of the length of the longest increasing
  subsequence of random permutations.
\newblock {\em J. Amer. Math. Soc.}, 12(4):1119--1178, 1999.

\bibitem[BDJ00]{BDJ99secondrow}
J.~Baik, P.~Deift, and K.~Johansson.
\newblock On the distribution of the length of the second row of a {Y}oung
  diagram under {P}lancherel measure.
\newblock {\em Geometric And Functional Analysis}, 10(4):702--731, 2000.

\bibitem[Bia01]{Bia01}
P.~Biane.
\newblock Approximate factorization and concentration for characters of
  symmetric groups.
\newblock {\em Internat. Math. Res. Notices}, 4:179--192, 2001.

\bibitem[Bil69]{Bil69}
P.~Billingsley.
\newblock {\em Convergence of Probability Measures}.
\newblock Wiley, 1969.

\bibitem[Bor99]{Borodin:limitJordan}
A.~M. Borodin.
\newblock The law of large numbers and the central limit theorem for the
  {J}ordan normal form of large triangular matrices over a finite field.
\newblock {\em J. Math. Sci.}, 96(5):3455--3471, 1999.

\bibitem[BOO00]{BOO00}
A.~Borodin, A.~Okounkov, and G.~Olshanski.
\newblock Asymptotics of {P}lancherel measures for symmetric groups.
\newblock {\em J. Amer. Math. Soc.}, 13(3):481--515, 2000.

\bibitem[Bri69]{Bri69}
D.~Brillinger.
\newblock The calculation of cumulants via conditioning.
\newblock {\em Ann. Inst. Statist. Math.}, 21:375--390, 1969.

\bibitem[C{\'S}09]{CollinsSniady09}
B.~Collins and P.~{\'S}niady.
\newblock Asymptotic fluctuations of representations of the unitary groups.
\newblock arXiv:0911.5546, 2009.

\bibitem[Gec93]{Geck93}
M.~Geck.
\newblock A note on harish-chandra induction.
\newblock {\em Manuscripta Math.}, 80:393--401, 1993.

\bibitem[GP00]{GP00}
M.~Geck and G.~Pfeiffer.
\newblock {\em Characters of Finite Coxeter Groups and Iwahori-Hecke Algebras},
  volume~21 of {\em London Mathematical Society Monographs}.
\newblock Oxford University Press, 2000.

\bibitem[HL80]{HL80}
R.~B. Howlett and G.~I. Lehrer.
\newblock Induced cuspidal representations and generalised {H}ecke rings.
\newblock {\em Invent. Math}, 58:37--64, 1980.

\bibitem[IO02]{IO02}
V.~Ivanov and G.~Olshanski.
\newblock Kerov's central limit theorem for the {P}lancherel measure on {Y}oung
  diagrams.
\newblock In {\em Symmetric Functions 2001: Surveys of Developments and
  Perspectives}, volume~74 of {\em NATO Science Series II. Mathematics, Physics
  and Chemistry}, pages 93--151, 2002.

\bibitem[Iwa64]{Iwa64}
N.~Iwahori.
\newblock On the structure of the {H}ecke ring of a {C}hevalley group over a
  finite field.
\newblock {\em J. Faculty Science Tokyo University}, 10:215--236, 1964.

\bibitem[Joh01]{JohanssonPlancherel}
K.~Johansson.
\newblock Discrete orthogonal polynomial ensembles and the {P}lancherel
  measure.
\newblock {\em Ann. of Math.}, 153(2):259--296, 2001.

\bibitem[Ker92]{KerovQHookWalk}
S.~V. Kerov.
\newblock $q$-analogue of the hook walk algorithm and random {Y}oung tableaux.
\newblock {\em Funct. Anal. Appl.}, 26(3):179--187, 1992.

\bibitem[Ker96]{Kerov96}
S.~V. Kerov.
\newblock A differential model for the growth of {Y}oung diagrams.
\newblock {\em Proceedings of the St. Petersburg Mathematical Society},
  4:165--192, 1996.
\newblock English translation: American Mathematical Society Translations,
  Series 2, 188, 1998.

\bibitem[Ker03]{Ker03}
S.~V. Kerov.
\newblock {\em Asymptotic representation theory of the symmetric group and its
  applications in analysis}, volume 219 of {\em Trans. Math. Monographs}.
\newblock AMS, 2003.

\bibitem[KOO98]{KOO97}
S.~V. Kerov, A.~Okounkov, and G.~Olshanski.
\newblock {The boundary of the Young graph with Jack edge multiplicities}.
\newblock {\em International Mathematics Research Notices}, 1998(4):173, 1998.

\bibitem[KO94]{KO94}
S.~V. Kerov and G.~Olshanski.
\newblock Polynomial functions on the set of {Y}oung diagrams.
\newblock {\em Comptes Rendus Acad. Sci. Paris S\'erie I}, 319:121--126, 1994.

\bibitem[KOV04]{KOV04}
S.~V. Kerov, G.~Olshanski, and A.~M. Vershik.
\newblock Harmonic analysis on the infinite symmetric group.
\newblock {\em Invent. Math.}, 158:551--642, 2004.

\bibitem[KV77]{KV77}
S.~V. Kerov and A.~M. Vershik.
\newblock Asymptotics of the {P}lancherel measure of the symmetric group and
  the limiting form of {Y}oung tableaux.
\newblock {\em Doklady AN SSSR}, 233(6):1024--1027, 1977.
\newblock English translation: Soviet Mathematics Doklady 18:527-531, 1977.

\bibitem[KV81]{KV81}
S.~V. Kerov and A.~M. Vershik.
\newblock Asymptotic theory of characters of the symmetric group.
\newblock {\em Funkts. Anal. i Prilozhen}, 15:15--27, 1981.
\newblock English translation: Funct. Anal. Appl. 15:246-255, 1982.

\bibitem[LS59]{LS59}
V.~P. Leonov and A.~N. Sirjaev.
\newblock On a method of semi-invariants.
\newblock {\em Theor. Prob. Appl.}, 4:319--329, 1959.

\bibitem[LS77]{LS77}
B.~F. Logan and L.~A. Shepp.
\newblock A variational problem for random {Y}oung tableaux.
\newblock {\em Adv. Math.}, 26:206--222, 1977.

\bibitem[Lus84]{Lus84}
G.~Lusztig.
\newblock {\em Characters of Reductive Groups over a Finite Field}, volume 107
  of {\em Ann. of Math. Studies}.
\newblock Princeton University Press, 1984.

\bibitem[Mac95]{Mac95}
I.~G. Macdonald.
\newblock {\em Symmetric functions and {H}all polynomials}.
\newblock Oxford Mathematical Monographs. Oxford University Press, 2nd edition,
  1995.

\bibitem[Mat99]{Mat99}
A.~Mathas.
\newblock {\em {I}wahori-{H}ecke algebras and {S}chur algebras of the symmetric
  group}, volume~15 of {\em University Lecture Series}.
\newblock Amer. Math. Soc., 1999.

\bibitem[M{\'e}l10]{Mel10}
P.-L. M{\'e}liot.
\newblock Gaussian concentration of the $q$-characters of the {H}ecke algebras
  of type {A}.
\newblock {I}n preparation, 2010.

\bibitem[Oko00]{Okounkov99}
A.~Okounkov.
\newblock Random matrices and random permutations.
\newblock {\em Internat. Math. Res. Notices}, 20:1043--1095, 2000.

\bibitem[Ram91]{Ram91}
A.~Ram.
\newblock A {F}robenius formula for the characters of the {H}ecke algebras.
\newblock {\em Invent. Math.}, 106:461--488, 1991.

\bibitem[RR97]{RR97}
A.~Ram and J.~B. Remmel.
\newblock Applications of the {F}robenius formulas for the characters of the
  symmetric group and the {H}ecke algebra of type {A}.
\newblock {\em Algebraic Combinatorics}, 5:59--87, 1997.

\bibitem[RRW96]{RRW96}
A.~Ram, J.~B. Remmel, and T.~Whitehead.
\newblock Combinatorics of the q-basis of symmetric functions.
\newblock {\em Journal of Combinatorial Theory}, 76:231--271, 1996.

\bibitem[SCc08]{Sage-Combinat}
The {S}age-{C}ombinat community.
\newblock {S}age-{C}ombinat: enhancing sage as a toolbox for computer
  exploration in algebraic combinatorics, 2008.
\newblock {\tt http://combinat.sagemath.org}.

\bibitem[{\'S}ni06]{Sniady06}
P.~{\'S}niady.
\newblock Gaussian fluctuations of characters of symmetric groups and of
  {Y}oung diagrams.
\newblock {\em Probab. Theory Relat. Fields}, 136:263--297, 2006.

\bibitem[S{\etalchar{+}}09]{sage}
W.\thinspace{}A. Stein et~al.
\newblock {\em {S}age {M}athematics {S}oftware ({V}ersion 4.2+)}.
\newblock The Sage Development Team, 2009.
\newblock {\tt http://www.sagemath.org}.

\bibitem[Str08]{Stra08}
E.~Strahov.
\newblock A differential model for the deformation of the {P}lancherel growth
  process.
\newblock {\em Adv. Math}, 217(6):2625--2663, 2008.

\bibitem[Tho64]{Tho64}
E.~Thoma.
\newblock Die unzerlegbaren, positive-definiten {K}lassenfunktionen der
  abz\"ahlbar unendlichen symmetrischen {G}ruppe.
\newblock {\em Math. Zeitschrift}, 85:40--61, 1964.

\end{thebibliography}

\end{document}